\documentclass[twoside]{article} 

\usepackage{amsmath,theorem,amssymb}

\newcommand{\R}{\mathbb{R}}
\newcommand{\C}{\mathbb{C}}
\newcommand{\Z}{\mathbb{Z}}
\newcommand{\N}{\mathbb{N}}

\newcommand{\A}{\textnormal{A}}
\newcommand{\B}{\textnormal{B}}
\newcommand{\rC}{\textnormal{C}}
\newcommand{\D}{\textnormal{D}}
\newcommand{\E}{\textnormal{E}}
\newcommand{\F}{\textnormal{F}}
\newcommand{\G}{\textnormal{G}}

\newcommand{\g}{{\mathfrak{g}}}    
\newcommand{\h}{{\mathfrak{h}}}    

\renewcommand{\L}{{\cal{L}}}    
\renewcommand{\O}{{\cal{O}}}    
\renewcommand{\P}{{\cal{P}}}    

\newcommand{\ord}{\textnormal{ord}}
\newcommand{\Ad}{\textnormal{Ad}}

\newcommand{\Proof}{\textsc{Proof:}\quad} 
\let\Box\undefined\newsymbol\Box 1203
\newcommand{\Endproof}{\hfill  $\Box$ \bigskip\par}
\newcommand{\Endproofshort}{\hfill  $\Box$ \par}

\newtheorem{Theorem}{Theorem}[section]
\newtheorem{Lemma}[Theorem]{Lemma}

\newtheorem{Cor}[Theorem]{Corollary}
\newtheorem{Prop}[Theorem]{Proposition}
{\theorembodyfont{\normalfont} \newtheorem{Rem}[Theorem]{Remark}}

\title{Weyl's character formula for non-connected Lie groups 
and orbital theory for twisted affine Lie algebras}

\author{Robert Wendt}

\begin{document}

\maketitle

\abstract{We generalize I.~Frenkel's orbital theory for 
non twisted affine Lie algebras
to the case of twisted affine Lie al\-ge\-bras using a character 
formula for certain non-connected compact
Lie groups.}

\section{Introduction}
In \cite{Frenkel}, Frenkel develops a generalization of Kirillov's
orbit theory for finite dimensional Lie groups \cite{Kirillov} in the
case of untwisted affine Lie algebras. In particular, he classifies
affine (co-)adjoint orbits of the underlying loop groups in terms of
conjugacy classes of connected compact Lie groups. Using this classification
and the theory of Wiener Integration on compact Lie  groups, Frenkel obtains
an interpretation of the character of a highest weight representation as an
integral over an associated orbit in the coadjoint representation.
\par
The aim of this paper is to generalize Frenkel's orbital theory to the case
of twisted affine Lie algebras.  Kleinfeld (\cite{Kleinfeld}) already made
a first step towards the adaptation of Frenkel's theory to this more
general  case, but he was unable to interpret the character formula as
an orbital integral.  In order to obtain an interpretation of the formulas
as orbital integrals we shall introduce certain non-connected 
compact Lie groups,
so called  principal extensions, into the geometrical picture. It turns out
that the  affine orbits of the adjoint representation of a twisted loop group
are parametrized by the conjugacy classes in the "outer component" of such a
group. With this result and a character  formula for non-connected Lie
groups at hand we are able to translate Frenkel's program to the twisted case.
I now give a brief description of the contents of this paper.
\par
\S \ref{integration} contains some facts about non-connected compact 
Lie groups,
some of them well known.    After constructing the principal extension
$\tilde G$ of a semisimple compact group $G$, we derive an analogue of 
the  Weyl
character formula for the connected components of $\tilde G$ not containing the
identity. 
The proof of this formula involves an analogue of 
the Weyl integral formula for $G$-invariant functions on $\tilde G$, 
and the character itself is
governed by the dual of a certain ``folded'' root system.
That is, if $R$
is the root system of $G$ and $\tau$ is a diagram  automorphism of $R$ 
then $\tau$
acts as an outer automorphism on $G$, and the characters on the connected
components of $\tilde G$ are governed by the root system $R^{\tau\vee}$, 
the dual to the
"folded" root system $R^\tau$ of the fixed point group $G^\tau$. 
It is worthwhile to
note that no group belonging to the root system $R^{\tau\vee}$ can, 
in general at
least, not be realized as a subgroup of $G$ in contrast to $G^\tau$.
The group $G$ acting on $\tilde G$ by conjugation, we shall view each
connected component of $\tilde G$ as a $G$-manifold.  
As another direct application
of the integral formula, we compute the radial component  
of the Laplacian on the
connected components of $\tilde G$ with respect to that $G$-action.
\par
In \S \ref{orbits}, we study affine orbits of the adjoint representation of
a twisted loop group $\L(G,\tau)$. By a slight alteration  of 
Frenkel's original
methods we see that for $G$ compact, every such orbit contains a 
constant loop and the
orbits in certain affine "shells" are parametrized by the $G$-orbits in the
connected component of $\tilde G$ containing $\tau$. To be more 
precise and using different
terminology, Frenkel regards a loop into the Lie  algebra 
$\g$ as a connection on a
principal, trivial fibre bundle over the circle $S^1$ with structure 
group $G$ and he
associates to this loop the monodromy of this connection which is an 
element of $G$.
The action of $\L(G)$ on an affine shell in the affine Lie algebra is then
given by gauge transformations and it is compatible with the 
$G$-conjugation on the
monodromies in $G$.  That way he obtains a well defined bijection 
from the adjoint
orbits of $\L(G)$ (of some fixed affine shell) to the set of conjugacy 
classes in $G$.
In the case of a twisted loop group $\L(G,\tau)$ where $\tau$ is a diagram
automorphism of $G$ the corresponding map will neither be surjective 
nor injective (cf.
\cite{Kleinfeld}). In this case, it is appropriate to replace the 
monodromy with the
``$\frac1r$-th monodromy'' (i.e. "monodromy" after $\frac1r$-th 
of the full circle) multiplied by $\tau$.
Here $r$ is the order of $\tau$. The gauge action of 
$\L(G,\tau)$ is then
compatible with conjugation in the component 
$G\tau$,  and we are able to classify the
affine adjoint orbits of $\L(G,\tau)$ in terms 
of conjugacy classes in $G\tau$.
\par
In \S\ref{character}, we use the theory of Wiener integration on a 
non-connected compact Lie group to rewrite the irreducible 
highest weight characters of a
twisted affine Lie algebra as an integral over the space 
of paths inside the connected
component of $\tilde G$ containing the element $\tau$. After a brief
summary of some results about affine Lie algebras and their representations in
\S\ref{affcharacters}, we shall show in \S\ref{poisson} 
and in \S\ref{heat} how the irreducible
highest weight characters of a twisted affine Lie algebra 
are linked to the fundamental
solution of the heat equation on a non-connected compact 
Lie group. At this point, the
characters on the connected component $G\tau$ enter the 
picture.  The fundamental
solution of the heat equation is used in 
\S\ref{wienermeasures} to define the Wiener
measure on the space of  paths in the connected 
component of a compact Lie group.
Computing a certain integral with respect to this 
measure, we can rewrite the affine characters
as an integral over a path space. In \S\ref{interpretation} 
we then  show, adopting the
original procedure of \cite{Frenkel}, how this integral 
can be interpreted as an integral over
a coadjoint orbit of the corresponding twisted loop 
group  thus  completing Frenkel's
program for twisted affine algebras as well.

\section{Integration and character formulas 
for non-connected Lie groups}\label{integration}

\subsection{Principal extensions and conjugacy classes}\label{extension}

Let $G$ be a simply connected semisimple compact Lie group, $T$ a maximal
torus in $G$,
and $R$ the root system of $G$ with respect to $T$.
If $W$ is the Weyl group of $R$ (and $G$) and $\Gamma$ is the group of diagram
automorphisms of the Dynkin diagram of $G$ then we have
$Aut(R)=W\rtimes\Gamma$.
Every $\tau\in\Gamma$ can be lifted to an automorphism of $G$ in the
following way.
Let $\Pi$ be a basis of $R$. If $\g$ is the Lie algebra of $G$ and $\g_\C$ its
complexification, we choose a set
$\{e_\alpha,f_\alpha,h_\alpha\}_{\alpha\in\Pi}$
of Chevalley generators of $\g_\C$ and set $\tau(e_\alpha)=e_{\tau(\alpha)}$.
This extends to a Lie algebra automorphism of $\g_\C$ which leaves  $\g$
invariant.
Since we chose $G$ to be simply  connected, $\tau$ can again be lifted to an
automorphism of $G$ leaving $T$ invariant. Thus we have defined a homomorphism
$\varphi:\Gamma\to Aut(G)$, and we can set $\tilde
G=G\rtimes_{\varphi}\Gamma$,
and call $\tilde G$ the principal extension of $G$. Obviously we have
$\tilde G/G=\Gamma$.
\par
If $G$ is not simply connected we have $G=\bar G/K$, where $\bar G$ denotes
the
universal covering group of $G$, and $K$ is some subgroup of the center of
$\bar G$.
Let $\Gamma_K$ be the subgroup of $\Gamma$ which leaves $K$ fixed. Then the
group
$G\rtimes_{\varphi}\Gamma_K$ is called the principal extension of $G$. The
principal
extensions of compact Lie groups are, in general, central extensions of the
automorphism groups of  the compact groups and play a crucial role in the
structure
theory of the non-connected  compact groups (cf. \cite{deSiebenthal}).
\par
$G$ acts on the components of $\tilde G$ by conjugation. The orbits of this
action
on the component $G\tau$, the connected component of $\tilde G$ containing
$\tau$,
are parametrized by a component of the space $S/W(S)$, where $S$ is a
Cartan  subgroup of
$\tilde G$ in the sense of \cite{BtD} such that $\tau\in S$. That is in our
cases
$S=T^\tau\times\langle\tau\rangle$ where $T^\tau$ is the connected
component containing
the identity of the $\tau$-invariant part of the maximal torus $T$. The
group $W(S)=N(S)/S$
is a finite group and is called the Weyl group belonging to $S$.
In particular, if $S_0=T^\tau$ is the connected component of $S$ containing
$e$, then
every element of $G\tau$ is conjugate under $G$ to an element of $S_0\tau$,
and two
elements of $S_0\tau$ are conjugate under $G$ if and only  if they are
conjugate in $N(S)$.
Furthermore, $S_0$ is regular in $G$, i.e. there is a unique
$\tau$-invariant maximal
torus of $G$ containing $S_0$, cf. e.g. \cite{BtD}.

\subsection{A 'Weyl integral formula' for $G\tau$}\label{integrationformula}

Let $G$ be a connected semisimple Lie group of type $\A_n$, $\D_n$, or
$\E_6$ and
$\tilde G$ its principal extension. Since the other connected Dynkin
diagrams do
not admit any diagram automorphisms, the principal extensions  of the
corresponding
compact Lie groups are trivial. We want to derive an analogue of the Weyl
integral
formula for $G$-invariant functions on the component $G\tau$  of $\tilde G$.
Let $S\subset \tilde G$ be a Cartan subgroup of $G$ containing $\tau$ such
that
$S/S_0$ is generated by $S_0\tau$, and let $(.,.)$ be the negative of the
Killing
form on $\g$. Since $G$ is semisimple, this gives an $\Ad(G)$ invariant scalar
product on $\g$, and $\g$ decomposes into a direct sum $\g=LT\oplus
L(G/T)$, where
$LT$ is the Lie algebra of $T$ and $L(G/T)$ its orthogonal complement. In
the same
way, $LT$ decomposes into $LT=LS_0\oplus L(T/S_0)$, hence we have
$\g=LS_0\oplus
L(T/S_0)\oplus L(G/T)$.
\par
There are normalized left invariant volume forms $d(gS_0)$, $ds$, and $dg$
on $G/S_0$,
$S_0$ and $G$ which are unique up to sign (i.e. orientation). The form $dg$
defines
a volume form on $G\tau$ by right translation with $\tau$, which will be
called $dg$
as well. Hence, on $G$, we have $|\Gamma_K|\cdot d\tilde g=dg$, where
$d\tilde g$ is the
normalized left invariant volume form on $\tilde G$.
\par
The projection $\pi:G\to G/S_0$ induces a map $D\pi:\g\to T_{eS_0}G/S_0$
which maps
$L(G/S_0):=L(G/T)\oplus L(T/S_0)$ isomorphically to the tangent space of
$G/S_0$ at
the point $eS_0$. Hence we can identify these spaces via $D\pi$. Let
$n=\dim G$ and
$k=\dim S_0$. Then $\pi^*d(gS_0)$ is a left invariant $(n-k)$-form on $G$.
Using the
$k$-form $ds_e\in\textnormal{Alt}^k(LS_0)$ we get
$\textnormal{pr}_2^*ds_e\in \textnormal{Alt}^k\g$
where $\textnormal{pr}_2:\g=L(G/S_0)\oplus LS_0\to LS_0$ is the second
projection.
The form $\textnormal{pr}_2^*ds_e$ defines a left-invariant $k$-form
$\beta$ on $G$ by
left translation, so  $\pi^*d(gS_0)\wedge\beta$ is a volume form on $G$.
Hence we have
$\pi^*d(gS_0)\wedge\beta=cdg$. We may chose the signs so that $c>0$, and it
is not hard
to see that in this case $c=1$ (cf. \cite{BtD}).
\par
There is a volume form
$\alpha=\textnormal{pr}_1^*d(gS_0)\wedge\textnormal{pr}_2^*ds$ on
$G/S_0\times S_0$. Identifying $\g$ with $L(G/S_0)\oplus LS_0$ and
evaluating the forms
at the unit element, one finds $\alpha_{(eS_0,e)}=dg_{(e)}$.
\par
Since every element of $G\tau$ is conjugate under $G$ to an element of
$S_0\tau$, the map
$$q:G/S_0\times S_0\to  G\tau,$$
$$(gS_0,s)\mapsto gs\tau g^{-1}$$
is surjective, and by the above we get $q^*dg=\det(q)\alpha$.

\begin{Lemma}
The functional determinant of the conjugation map $q$ is given by
$$\det(q)(gS_0,s)=\det(\Ad|_{L(G/S_0)}(s\tau)^{-1}-I|_{L(G/S_0)}).$$
\end{Lemma}
\Proof Similar to the proof of Prop. IV, 1.8 in \cite{BtD}.
\footnote{For more details on the proofs in this section see
\cite{Wendt}}
\Endproof

\begin{Lemma}
Let $z$ be a generator of $S$. Then
\renewcommand{\theenumi}{\roman{enumi}}
\renewcommand{\labelenumi}{\textnormal{(\theenumi)}}
\begin{enumerate}
\item
$|q^{-1}(z)|=|W(S)|$.
\item
For $(g_1S_0,s_1),(g_2S_0,s_2)\in q^{-1}(z)$ we have
$$\det(q)(g_1S_0,s_1)=\det(q)(g_2S_0,s_2).$$
\item
There exists a generator $z$ of $S$ such that $q$ is regular in each \\
 $(gS_0,s)\in q^{-1}(z)$.
\end{enumerate}
\end{Lemma}
\Proof Similar to the proof of Prop. IV, 1.9 in \cite{BtD}.
\Endproof

From this we obtain the mapping degree of $q$,
$\textnormal{deg}(q)=\textnormal{sign}(\textnormal{det}(q))\cdot|W(S)|$.
Thus using Fubini's theorem, one gets

\begin{Prop}
Let $f:G\tau\to\R$ be a $G$-invariant function.Then
$$\int_{G\tau}f(g)dg=\frac{1}{|W(S)|}\cdot\int_{S_0}f(s\tau)
\cdot|\det(\Ad|_{L(G/S_0)}(s\tau)^{-1}-I|_{L(G/S_0)})|ds.$$
\end{Prop}

In order to obtain an analogue of the classical Weyl integral formula, one has
to calculate the functional determinant of $q$ in  terms of the root system
of $G$.
We adapt the notation used in \cite{BtD}. That is, if $\alpha\in LT^*$ is an
infinitesimal root of $G$, then $\vartheta_\alpha$ denotes the
corresponding global
root $T\to S^1$. So for $H\in LT$ one has $\vartheta_\alpha\circ\exp(H)=
e^{2\pi i \alpha(H)}$.
Setting $e(x)=e^{2\pi i x}$, we get $\vartheta_\alpha=e(\alpha)$.
\par
Now we consider the action of $(\Ad(s\tau)-I)$ on $L(G/S_0)_\C=\C\otimes
L(G/S_0)$.
As before, $L(G/S_0)_\C$ decomposes into two orthogonal subspaces
$L(G/S_0)_\C=L(T/S_0)_\C\oplus
L(G/T)_\C$. Let $X\in L(T/S_0)$ be an eigenvector of $\Ad(\tau)$. A short
calculation shows
$$(\Ad(s\tau)-I)(X)=-X+\gamma X,$$ where $\gamma=\pm 1$ if $\tau^2=e$, and
$\gamma$ is a
third root of unity if $\tau^3=e$. We can choose $(gS_0,s)$ to be a regular
point of $q$,
so $\gamma\neq 1$ and $ (\Ad(s\tau)-I)(X)=-2X$ if $t^2=e$, and $
(\Ad(s\tau)-I)(X)=(\gamma-1)X$
if $\tau^3=e$.
\par
Now $L(G/T)_\C$ decomposes into the direct sum of root spaces
$$ L(G/T)_\C = \bigoplus_{\alpha\in R}L_\alpha,$$
and for $X\in L_\alpha$, $s\in S_0$, one has
$$\Ad(s)(X)=\vartheta_\alpha(s)(X).$$
\par
Let $\bar\alpha=\alpha|_{LS_0}$ for $\alpha\in R$. It is a well known fact
that the
set $R^\tau=\{\bar\alpha|\alpha\in R\}$ is a (not necessarily reduced) root
system.
The relation between the type of $R$ and the type of $R^\tau$ is shown in
the  following
table:

$$
\begin{array}{c|c|c|c|c|c|}
R & \A_{2n-1} & \A_{2n} & \D_n (n\geq4) & \D_4 & E_6 \\
\hline
\textnormal{ord}(\tau) & 2 & 2 & 2 & 3 & 2  \\
\hline
R^\tau & \rC_n & \B\rC_n & \B_{n-1} & \G_2 & \F_4
\end{array}
$$

So if $R$ is of type $\A_{2n-1}$, $ \D_n$, or $E_6$, then $R^\tau$ is a
reduced root system.
It is easy to see that in this case, $\bar\alpha$ is a long root of
$R^\tau$ if and only if
$\tau(\alpha)=\alpha$. Otherwise $\bar\alpha$ is a short root of $R^\tau$.
If $R$ is of type
$\A_{2n}$, then the root system $R^\tau$ is not reduced and three distinct
root lengths occur.
In this case, $\bar\alpha$ is a long root in $R^\tau$ if $\alpha$ is
invariant under $\tau$.
If $\alpha$ and $\tau(\alpha)$ are orthogonal to each other, then
$\bar\alpha$ is a root of
medium length in $R^\tau$, and otherwise $\bar\alpha$ is a short root in
$R^\tau$.  (Remember
that the root system $\B\rC_n$ is the union of two root systems of types
$B_n$ and $C_n$
such that the long roots of $B_n$ coincide with  the short roots of $C_n$.)
\par
Now let $G$ be of type $\A_{2n-1}$, $ \D_n$, or $E_6$, and consider the
case $\tau^2=e$.
We have seen that $\tau$ defines a Lie algebra automorphism via
$\tau(X_\alpha)=X_{\tau(\alpha)}$
for $\alpha\in\Pi$ and $X_\alpha\in L_\alpha$.  Extending this to the
entire Lie algebra, one
gets $\tau(X_\alpha)=X_{\tau(\alpha)}$ for all $\alpha\in R$. The
eigenvectors of
$\Ad(\tau)$ are the following: If $\alpha$ is invariant under $\tau$ then
$X_\alpha$ is an
eigenvector with eigenvalue 1. If $\alpha$ is not invariant under $\tau$,
there are two
eigenvectors $X_\alpha\pm X_{\tau(\alpha)}$ of eigenvalue $\pm1$. Thus for
$\tau(\alpha)=\alpha$
and $X=X_\alpha$ we have
$$(\Ad(s\tau)^{-1}-I)(X)=(\vartheta_{\bar\alpha}(s^{-1})-1)X,$$
and
$$(\Ad(s\tau)^{-1}-I)(X)=(\pm\vartheta_{\bar\alpha}(s^{-1})-1)X$$
for $\tau(\alpha)\neq\alpha$ and $X=X_\alpha\pm X_{\tau(\alpha)}$
respectively. This yields
\begin{multline*}
\det(\Ad|_{L(G/S_0)}(s\tau)^{-1}-  I|_{L(G/S_0)}) = \\
 (-2)^{\dim(T/S_0)}\cdot\prod_{\underset{\tau(\alpha)=\alpha}{\bar\alpha\in R^\tau}}
(\vartheta_{\bar\alpha}(s^{-1})-1)\\
 \cdot\prod_{\underset{\tau(\alpha)\neq\alpha}{\bar\alpha\in
R^\tau}}(\vartheta_{\bar\alpha}(s^{-1})-1)
(-\vartheta_{\bar\alpha}(s^{-1})-1)
\end{multline*}
Multiplying each factor by
$-1$, using the remark above on the relative length of the 
$\bar\alpha$ as well as the equality
$$(1-\vartheta_{\bar\alpha})(1+\vartheta_{\bar\alpha})=(1-\vartheta_{2\bar\alpha
}),$$
this becomes
\begin{align*}
  & = (-2)^{\dim(T/S_0)}\cdot\prod_{\underset{\bar\alpha\text{
long}}{\bar\alpha\in R^\tau}}
(1-\vartheta_{\bar\alpha}(s^{-1}))
        \prod_{\underset{\bar\alpha\text {short}}{\bar\alpha\in
R^\tau}}(1-\vartheta_{\bar2\alpha}(s^{-1}))    \\
  & =  (-2)^{\dim(T/S_0)}\cdot\Delta(s^{-1})\bar\Delta(s^{-1}),
\end{align*}
with
$$\Delta(s)=\prod_{\bar\alpha\in
R^{\tau\vee}_+}(1-\vartheta_{\bar\alpha}(s)).$$
Here $R^{\tau\vee}$ denotes the dual root system of $R^\tau$ which is given by
$\bar\alpha^\vee=\frac{2\bar\alpha}{\langle\bar\alpha,\bar\alpha\rangle}$ for
$\bar\alpha\in R^\tau$ and $\langle.,.\rangle$ is a multiple of the
Killing form such that $\langle\bar\alpha,\bar\alpha\rangle=2$ for a long root
$\bar\alpha\in R^\tau$.
\par
\medskip
If $G$ is of type $\D_4$ and $\tau^3=e$, a similar calculation gives
$$\det(\Ad|_{L(G/S_0)}(s\tau)^{-1}-  I|_{L(G/S_0)}) = 3\cdot
\Delta(s^{-1})\bar\Delta(s^{-1}),$$
with $\Delta(s)$ as above. Observe that in this case $R^\tau$ is of type
$\G_2$, so $\dim T/S_0=2$.
\par
If $G$ is of type $A_{2n}$ we have to be more careful since  $R^\tau$
is not reduced and
three different root lengths occur. Also, in this case 
the Lie algebra automorphism
$\tau$ is slightly more
complicated. For $X_\alpha\in L_\alpha$ we have
$$\tau(X_\alpha)=(-1)^{1+\textnormal{ht}(\alpha)}X_{\tau(\alpha)}.$$
Now $\tau(\alpha)=\alpha$ implies that $\textnormal{ht}(\alpha)$ is even
and a similar
calculation yields
\begin{multline*}
\det(\Ad|_{L(G/S_0)}(s\tau)^{-1}-  I|_{L(G/S_0)}) = \\
  (-2)^{\dim(T/S_0)}\cdot\prod_{\underset{\bar\alpha\text{ long}}{\bar\alpha\in
R^\tau}}(1+\vartheta_{\bar\alpha}(s^{-1}))
  \cdot\prod_{\underset{\bar\alpha\text{ middle}}{\bar\alpha\in
R^\tau}}(1-\vartheta_{2\bar\alpha}(s^{-1}))\\
   \cdot\prod_{\underset{\bar\alpha\text{ short}}{\bar\alpha\in
R^\tau}}(1-\vartheta_{2\bar\alpha}(s^{-1}))
\end{multline*}
But the length of the long roots in $\B\rC_n$ is twice the length of the
short roots. So
we can put these together to obtain $$\det(\Ad|_{L(G/S_0)}(s\tau)^{-1}-
I|_{L(G/S_0)}) =
(-2)^{\dim(T/S_0)}\cdot\Delta(s^{-1})\bar\Delta(s^{-1}),$$ with
$$\Delta(s)=\prod_{\bar\alpha\in R^{1}_+}(1-\vartheta_{\bar\alpha}(s)).$$
Here $R^1=\{2\bar\alpha|\bar\alpha\in\B\rC_n,\bar\alpha\text{ long}\}\cup
\{2\bar\alpha|\bar\alpha\in\B\rC_n,\bar\alpha\text{ middle}\}$. This is a
root system
of type $\rC_n$.
\par
\medskip
Before stating the integral formula for $G\tau$, we have to compare the
different Weyl groups
involved. Let $T$ be the maximal torus of $G$ such that $S_0\subset T$ and
let $W(T)=N_G(T)/T$
be the usual Weyl group of $G$. If we set $W^\tau=\{w\in W(T)|\tau
w\tau^{-1}=w\}$, then
$W^\tau$ is the Weyl group of the root system $R^\tau$ (and  also of course
of its dual
$R^{\tau\vee}$).
\begin{Prop}\label{weylgruppen}
Let $W^\tau$ be as above. Then
\renewcommand{\theenumi}{\roman{enumi}}
\renewcommand{\labelenumi}{\textnormal{(\theenumi)}}
\begin{enumerate}
\item
There exists a split exact sequence
$$e\to(T/S_0)^\tau\to W(S)\to W^\tau \to e.$$
Here $(T/S_0)^\tau$ denotes the fixed point set under conjugation with $\tau$.
\item
We have
$$|W^\tau|=
\begin{cases}
\frac1{2^{\dim(T/S_0)}}\cdot|W(S)| \quad & \text{ if }\tau^2=e\\
\frac13\cdot|W(S)| & \text{ if }\tau^3=e
\end{cases}$$
\end{enumerate}
\end{Prop}
\Proof
Observe that in our cases we have $W(S)=(N(S)\cap G)/G$ (this is not true
for general
non-connected Lie groups), so one can define a map $\varphi:W(S)\to W(T)$,
$gS_0\mapsto gT$.
This map is well defined and one has $Im(\varphi)\subset W^\tau$. The rest
is done by a
calculation in the Lie algebra of $T$. For details see \cite{Wendt} (also
cf. \ref{weylgroupsII}).
To see that the sequence splits, observe that $S_0$ is by construction a
maximal torus in the
connected component of $G^\tau$ containing $e$,  and $W^\tau$ is the
corresponding Weyl group.
\Endproof
Putting everything together, one gets
\begin{Theorem}['Weyl integral formula' for $G\tau$]\label{Integralformel}
Let $f:G\tau\to\R$ be a function which is integrable and invariant under
conjugation by $G$. Then
$$\int_{G\tau}f(g)dg=\frac1{|W^\tau|}\int_{S_0}f(s\tau)\Delta(s)\bar\Delta(s)ds,
$$
with
$$\Delta(s)=\prod_{\bar\alpha\in R^1_+}(1-\vartheta_{\bar\alpha}(s)).$$
Here $R^1$ denotes the root system $R^{\tau\vee}$ if $R^\tau$ is reduced and
$R^1=\{2\bar\alpha|\bar\alpha\in\B\rC_n,\bar\alpha\text{ long}\}
\cup\{2\bar\alpha|\bar\alpha\in\B\rC_n,\bar\alpha\text{ middle}\}$ is
of type $\rC_n$ if $R$ is of type $A_{2n}$.
\end{Theorem}

\subsection{Applications of the integral formula}\label{applications}

Let $G$ be a compact connected semisimple Lie group of type $\A_n$, $D_n$,
or $E_6$ as above,
and let $\bar G\subset \tilde G$ be any non trivial subextension of $G$.
As a first application of the integral formula we compute the irreducible
characters of
$\bar G$ on the component $G\tau$ for $\tau\in\Gamma$. Let $T$ be the
$\tau$-invariant
maximal torus of $G$, and let $S$ be a Cartan subgroup of
$G$ such that $S/S_0$ is generated by $\tau S_0$. In the case $\D_4$, this
notation is
not quite unique since in general more than one diagram automorphism occurs,
but it will
always be clear which respective Cartan subgroup is being used at the moment.
\par
Since $S_0$ is regular
in $G$, we can choose a Weyl chamber
$K\subset LT^*$ such that $K\cap LS_0^*$ is not empty and we let $\bar K$
denote its closure.
Then the set $K^\tau=K\cap LS_0^*$ is a Weyl chamber in $LS_0^*$ with
respect to the root system
$R^\tau$. Furthermore,  let $W=W(T)$ be the Weyl group of $G$, let $I$
denote the lattice
$I=ker(\exp)\cap LT$, and $I^*\subset LT^*$ its dual.
\par
For linear forms $\lambda\in LT^*$ and $\mu\in LS_0^*$ we define the
alternating sums
$$A(\lambda)=\sum_{w\in W(T)}\epsilon(w)\cdot e(w\lambda)$$
and
$$A^\tau(\mu)=\sum_{w\in W^\tau}\epsilon(w)\cdot e(w\mu)$$
respectively. Here we have set $\epsilon(w)=(-1)^{\textnormal{length}(w)}$.
Note that for $w\in W^\tau$ the two $\epsilon(w)$ in the equations above
do not necessarily coincide since they come from 
the presentations of $w$ as an 
element of two different Weyl groups.
In this notation, $A^\tau$ is a complex valued function on 
$LS_0$ which is alternating with respect to $W^\tau$.
\par
Now let $R^1$ be the root system used in Theorem
\ref{Integralformel}. As usual, we set
$$\rho=\frac12\sum_{\alpha\in R_+} \alpha,$$
and
$$\rho^\tau=\frac12\sum_{\bar\alpha\in R^1_+} \bar\alpha.$$
Now we can define functions $\delta:LT\to\C$, resp. $\delta^\tau:LS_0\to\C$ via
$$\delta=e(\rho)\cdot \prod_{\alpha\in R_+}(1-e(-\alpha))$$
and
$$\delta^\tau=e(\rho^\tau)\cdot \prod_{\bar\alpha\in
R^1_+}(1-e(-\bar\alpha)).$$
The function $\delta^\tau\cdot\bar\delta^\tau$ factorizes through $\exp$,
and we have
$\delta^\tau\cdot\bar\delta^\tau=(\Delta\cdot\bar\Delta)\circ\exp$, where
$\Delta$ is as in Theorem
\ref{Integralformel}. With this notation the classical Weyl character
formula for the irreducible
character $\chi_\lambda$ of $G$ belonging to the highest weight $\lambda\in
I^*\cap\bar K$ reads
$$\chi_\lambda=A(\lambda+\rho)/\delta.$$
\par
Now if $\mu\in I^*\cap LS_0^*$ then $A^\tau(\mu+\rho^\tau)/\delta^\tau$ can
be extended uniquely to
a function on $LS_0$. This function factors through $\exp$ (cf.
\cite{BtD}). In this way
$A^\tau(\mu+\rho^\tau)/\delta^\tau$ can be considered as a function on $S_0$.
\par
For an arbitrary character $\chi$  of $\bar G$  we define the function
$\chi^\tau:S_0\to\C$ via
$\chi^\tau(s)=\chi(s\tau)$. So on the connected component $G\tau$ of $\bar
G$, the character $\chi$
is determined by $\chi^\tau$. Now we can state the analogue of the Weyl
character formula for the
subextensions $\bar G\subset \tilde G$. First we consider the case when
$\bar G$ consists of two
connected components. So only one diagram automorphism $\tau$ is involved,
and $\tau^2=e$.

\begin{Theorem}\label{Charakterformel}
There exists an irreducible character $\tilde\chi_\lambda$ of $\bar G$ for
each
$\lambda\in I^*\cap \bar K$.
If $\lambda\not\in LS_0^*$, then
$$\tilde\chi_\lambda |_G = \chi_\lambda + \chi_{\tau(\lambda)},$$
and
$$\chi_\lambda |_{G\tau}\equiv 0.$$
Here $\chi_\lambda$ denotes the irreducible character of $G$ of highest
weight $\lambda$
\par
For each $\lambda\in LS_0^*$,  there exist two irreducible characters of
$\bar G$ associated to
$\lambda$, and we have
$$\tilde\chi_\lambda |_G = \chi_\lambda,$$
and
$$\tilde\chi^\tau_\lambda = \pm A^\tau(\mu+\rho^\tau)/\delta^\tau.$$
\end{Theorem}
If $\bar G$ consists of three connected components, the character formula
is essentially the same,
except that there are three irreducible characters for each $\lambda\in
LS_0^*$. If $G$ is of type
$\D_4$ and $\bar G=\tilde G$, then  one can use the character formulas
above together with some
information about the conjugacy classes in $S_3$ to determine all
irreducible characters. Since
this will not be needed in the sequel, we will omit the statement of the
result.
\par
\medskip
%
%
{\sc Remarks on the proof:}\quad
The proof of theorem \ref{Charakterformel} is essentially the same as the
proof of the classical
Weyl character formula in \cite{BtD}. The integral
formula is used along with the orthogonality relations for irreducible
characters to show that
the irreducible characters of $\bar G$ must have the given form. Then one
can apply the Peter-Weyl
theorem to see that each of the functions above must be an irreducible
character of $\bar G$.
For more details on this see \cite{Wendt}.
\Endproofshort
\begin{Rem} Kostant \cite{Kostant} states a character formula for 
non-connected complex 
Lie groups which is not quite as explicit as theorem \ref{Charakterformel}. 
In particular, 
he gives the character on the component $G\tau$ as a function on $T\tau$, 
where $T$ is a 
maximal torus in $G$. So the different root systems do not appear explicitly 
(although it is not hard to derive theorem  \ref{Charakterformel} 
from his formula).
\par
The formula itself was discovered before by Jantzen \cite{Jantzen}, who calculated 
the trace of the outer automorphism $\tau$ on the 
weight spaces of an irreducible representation with invariant highest weight of a 
semisimple algebraic group, and  
by Fuchs et al. \cite{Fuchs-Schweigert} who studied the characters of 
``$\tau$-twisted'' representations of a generalized Kac-Moody algebra.
\end{Rem}
As in the classical case, there is  a
Weyl denominator formula:
\begin{Cor}
With the same notation as above we have
$$\delta^\tau =  A^\tau(\rho^\tau).$$
\end{Cor}
\par
\begin{Rem}
It is interesting to note that the group belonging to the root system $R^1$
can, in general,
not be realized as a subgroup of $G$. For example,
$SO(2n+1)$, or its covering group $Spin(2n+1)$, which are the groups with
root system $\B_n$
can not be realized as subgroup of $SU(2n)$ which belongs to the root
system $\A_{2n-1}$.
\end{Rem}
\par
As a second application of theorem \ref{Integralformel} we can derive a
formula for the radial
part of the Laplacian on $G\tau$ with respect to the $G$-action by
conjugation. The negative of
the Killing form on
$\g_\C$ defines a positive definite scalar product on $\g$ which defines a
biinvariant Riemannian
metric on $G$ by left translation. We can use right multiplication by
$\tau^{-1}$ to pull back this Riemannian metric to $G\tau$. So $G$ and
$G\tau$ are isomorphic
as Riemannian manifolds. This metric is invariant under the $G$-action on
$G\tau$. Let
$\Delta_{G\tau}$ be the Laplacian on $G\tau$ with respect to this metric,
and let $\Delta_{S_0}$
be the Laplacian on $S_0$. Now we can use the general theory of radial
parts of invariant
differential operators
\cite{Helgason} Part III, Theorem 3.7 (also cf. Prop. 3.12)  
to get the following Proposition.

\begin{Prop}\label{radial}
Let $f:G\tau\to\R$ be a $G$-invariant function, and $f^\tau:G\to\R$ be
given by $f^\tau(g)=f(g\tau)$.
Let $\delta^\tau$ and $\rho^\tau$ be as above. Then we have
$$\delta^\tau\cdot\Delta_{G\tau}(f) =
(\Delta_{S_0}+\Vert\rho^\tau\Vert^2)(\delta^\tau\cdot f^\tau),$$
as functions restricted to $LS_0$. Here $\Vert.\Vert$  is the metric on
$LS_0^*$ induced by the
negative of the Killing form.
\end{Prop}

\section{The affine adjoint representation of a twisted loop
group}\label{orbits}
\subsection{Affine Lie algebras}\label{topologies}

Let $\g_\C$ be a complex semi simple Lie algebra, with compact involution
$\omega$
and let $\g$ be a compact form of $\g_\C$. That is,
$i\g=\{x\in\g_\C\enspace | \enspace\omega(x)=x\}$.
The loop algebra  $\L\g_\C$  (resp. $\L\g$) is the algebra of
$C^\infty$-maps from $S^1$
to $\g_\C$ (resp. $\g$). It is a Lie algebra under pointwise Lie bracket.
If $\g_\C$ is a complex simple Lie algebra of type $X_n$
and if the circle $S^1$ is parametrized by the real line $\R$ via the
exponential $e(t)=e^{2\pi i t}$
then the untwisted affine
Lie algebra\footnote{We shall adhere to some slight abuse of terminology,
here and in the sequel. The affine algebras in the sense of \cite{Kac}  are
based on algebraic loops,
with finite Fourier expansion. Our algebras may be viewed as completions of
these algebraic ones, cf.
below. The notation for types will be that of \cite{Kac}.}
$\tilde\g_\C$ of type $X_n^{(1)}$ is given by
$$\tilde\g_\C = \L\g_\C \oplus\C C\oplus \C D,$$
with Lie bracket
\begin{multline*}
[x(\cdot)+aC+bD,x_1(\cdot)+a_1C+b_1D]\\
= [x,x_1](\cdot) + bx_1'(\cdot)-b_1x'(\cdot)+(x'(\cdot),x_1(\cdot))C.
\end{multline*}
Here $x'(t)=\frac{dx(t)}{dt}$, $[x,x_1](t)=[x(t),x_1(t)]$, and
$$(x(\cdot),x_1(\cdot))=\int_0^1(x(t),x_1(t))dt,$$
where $(.,.)$ under the integral sign denotes the Killing form on $\g_\C$.
\par
An invariant bilinear form $(.,.)$ on $\tilde\g$ is given by
$$(x(\cdot)+aC+bD,x_1(\cdot)+a_1b+b_1D) = \int_0^1(x(t),x_1(t))dt+ab_1+a_1b.$$
\par
We obtain a so-called compact form
$\tilde \g$ of $\tilde\g_\C$ by considering
$\tilde\g = \L\g \oplus\R C\oplus \R D$.
\par
If $\tau$ is the outer automorphism of $\g_\C$ considered in \S\ref{extension}
and $\ord(\tau)=r$, then the twisted loop algebras $\L(\g_\C,\tau)$ and
$\L(\g,\tau)$
are given by $$\L(\g_{(\C)},\tau)=\{x\in\L
\g_{(\C)}|\tau(x(t))=x(t+1/r),\text {for all }t\in[0,1]\}.$$
Now if $\g_\C$ is a simple Lie algebra of type $X_n$ with $X_n=\A_n$,
$\D_n$ or $E_6$,
then the twisted affine Lie algebra $\tilde\g_\C$ of type $X_n^{(r)}$ is
given by
$$\tilde\g_\C = \L(\g_\C,\tau) \oplus\C C\oplus \C D$$
with the same Lie bracket as above.
The invariant bilinear form $(.,.)$ on the corresponding untwisted affine
Lie algebra
yields an invariant bilinear form  on the twisted affine Lie algebra
$\tilde\g_\C$ by
restriction and it is denoted by the same symbol.
The compact form $\tilde \g$ of a twisted affine Lie algebra is obtained in
the  same
way as in the untwisted case.
\par
We now define the $C^\infty$-topology on the Lie algebras $\tilde\g$. The
topology on
$\tilde\g_\C$ is obtained by viewing  it as a direct sum of two copies of
$\tilde\g$.
If $\L\g$ is an untwisted loop algebra, one defines the $C^\infty$-topology
on  $\L\g$
via the set of semi-norms
$$p_n(x)=\underset{\underset{t\in[0,1]}{s\leq n}}{\textnormal{sup}}\vert
\frac{d^sx(t)}{dt^s}\vert,
\quad n=1,2,\ldots,\quad x\in \L\g.$$
With respect to this topology $\L\g$ is complete. It extends immediately to
the enlarged
$\tilde\g$ as well as to the twisted subalgebras.We may view these
topological algebras as
$C^\infty$-completions of the algebraic loop algebras as e.g. in \cite{Kac}.
\par
%
%
%
%
\medskip
For later applications, we introduce a second topology on the spaces $\L(\g)$
and $\L(\g,\tau)$. The completions of these spaces with respect to the new
topology 
will not carry a Lie algebra structure any more (cf. \cite{Frenkel}). 
As before, let
$(.,.)$ denote the Killing form on $\g_\C$, and let $(.,.)_\g$
denote the negative of the Killing form on $\g$. Then 
$(x,y)_\g=\int_0^1(x(t),y(t))_\g dt$ gives a scalar product on
$\L(\g)$, and by restriction we get a scalar product on $\L(\g,\tau)$ 
as well. Let $\L(\g,\tau)(L_2)$ denote the completion of $\L(\g,\tau)$
with respect to the metric induced by the scalar product. 
After extending the scalar product to the completion, 
$\L(\g,\tau)(L_2)$ is a Hilbert space.
%
%
%
%

\subsection{Loop groups and their affine adjoint
representation}\label{adjointrep}

From now on, we will write $G_{(\C)}$ if we mean either a compact group
$G$ or the
corresponding complex group $G_\C$, and analogously we write
$\tilde\g_{(\C)}$ for
the associated Lie algebras.
\par
If $G_{(\C)}$  is a simply connected compact (complex) Lie group then the
corresponding
untwisted loop group $\L G_{(\C)}$  is defined to be the topological group
of $C^\infty$
mappings from $S^1$ to  $G_{(\C)}$  with pointwise multiplication and the usual
$C^\infty$-topology. If $\tau$ is one of the outer automorphisms of
$G_{(\C)}$
considered in \S\ref{extension} and $\ord(\tau)=r$ then the twisted loop group
$\L (G_{(\C)},\tau)$ is the subgroup
$$\L (G_{(\C)},\tau)=\{x\in \L G_{(\C)}| \tau(x(t))=x(t+1/r)\text{ for all
} t\in[0,1]\}.$$
\par
Then $\L \g_{(\C)}$ (resp. $\L (\g_{(\C)},\tau)$) may be viewed as the Lie
algebra of
$\L G_{(\C)}$ (resp. $\L (G_{(\C)},\tau)$) and there are natural adjoint
actions of these
groups on their Lie algebras (by pointwise finite-dimensional adjoint action).
However, for purposes of representation theory it is essential to consider
not this
action but the adjoint action of the affine Kac-Moody groups on the affine
Lie algebras.
These groups are given, similar to the Lie algebra case, as central
extensions of semidirect
products
$$e\to S^1\to\hat G \to \L G\rtimes S^1 \to e,$$
where the circle in the semidirect product acts on $\L G$ by rotation on
the loops (in
the complex case the circles $S^1$ are usually replaced by $\C^*$), and
accordingly in
the twisted cases. See \cite{PS} for a construction of  these groups.
Obviously the center
of $\hat G$ acts trivially in the adjoint representation, thus it is
sufficient to consider
the action of the quotient group $\L G\rtimes S^1 $. It will turn out that
the adjoint action of
the rotation group $S^1$ is without relevance for our purposes (i.e. in the
case that an
$\L G$-orbit contains a constant loop , see below, this action preserves
the $\L G$-orbit).
Therefore, it is sufficient to look at the adjoint action of only the loop
group $\L G$ on
$\tilde\g$. We call that action, which significantly differs from the usual
adjoint action
of $\L G$ on its Lie algebra $\L g$, the affine adjoint action and denote
it by $\widetilde{\Ad}$.
Exploiting the natural exponential mapping from $\L g$ to $\L G$ and
working with a fixed faithful
matrix representation of $G$, Frenkel was able to determine the exact form
of the affine adjoint
action in the  untwisted case (\cite{Frenkel}).
\begin{Prop}\label{ad}
Let $g\in\L G_{(\C)}$, $y\in \L\g_{(\C)}$ and $a,b\in\C$ (resp. $\R$), then
$$
\widetilde{\Ad}\enspace g(aC+bD+y) =\tilde aC+bD+gyg^{-1}-bg'g^{-1}
$$
with $\tilde a=a+(g^{-1}g',y)-\frac b2(g'g^{-1},g'g^{-1})$ and
$g'(t)=\frac{dg(t)}{dt}$.
\end{Prop}
The affine adjoint representation of a twisted loop group
$\L(G_{(\C)},\tau)$ on the
corresponding twisted affine Lie algebra $\tilde\g_{(\C)}$ is obtained by
restriction
of the adjoint representation of $\L G_{(\C)}$ on the corresponding
untwisted affine
Lie algebra $\tilde\g'_{(\C)}$ to $\L(G_{(\C)},\tau)$ and
$\tilde\g_{(\C)}$. Hence the
formula of Proposition \ref{ad} remains true for the twisted cases, as well.

\subsection{Classification of affine orbits}\label{orbitclassification}

In the case of an untwisted affine Lie algebra Frenkel and Segal have
classified affine
adjoint orbits of the loop group $\L G_{(\C)}$ on the affine Lie algebra
$\tilde\g_{(\C)}$
in terms of conjugacy classes of the group $G_{(\C)}$ (\cite{Frenkel},
\cite{PS}).
By a slight alteration of Frenkel's original methods, we can extend this
classification
to the twisted loop groups $\L (G_{(\C)},\tau)$.
Technically, the study of the affine adjoint action (on elements with
$b\neq0$) is
equivalent to that of transformations of ordinary differential equations
on the real
line  with periodic coefficients or, in more advanced terminology, to that
of the action
of the gauge group on connections in a trivial fibre bundle over the circle
$S^1$. In
the context of twisted loops we shall, in addition, have  to work with
differential
equations with 'twisted periodic' coefficients. To obtain a unified
statement of the
results, we shall allow the twisting automorphism $\tau$ to be the identity
on $G_{(\C)}$.
In this case, the group $\L (G_{(\C)},\tau)$ is just the untwisted group
$\L G_{(\C)}$.
\par
\medskip
Consider first the system of linear differential equations
$$z'(t)=z(t)x(t),$$
where $x(t),z(t)\in M_n(\C)$ for all $t\geq0$, and $x(t)$ is continuous in
$t$. A
fundamental result from the theory of differential equations secures the
existence
of a unique solution  $z(t)$ of the above equation with $z(0)=I_n$. This
solution
$z$ is usually called the fundamental solution of the differential equation.
\par
Now let $x$ be twisted periodic, that is, $x(t+1/r)=\tau x \tau^{-1}$ for some
invertible matrix $\tau$ with $\tau^r=I_n$ and all $t\geq 0$. If $z$ is the
fundamental solution of $z'=zx$, then obviously $z_1(t)=\tau^{-1}z(t+1/r)\tau$
is another solution.  Hence there exists a matrix $\widetilde M(x)$ such that
$z_1(t)=\widetilde M(x)z(t)$. Since we have chosen $z(0)=I_n$, we get
$\widetilde M(x)=z_1(0)=\tau^{-1}z(1/r)\tau$. Now let $M(x):=z(1/r)$ be the
"$\frac1r$-th monodromy" of the differential equation $z'=zx$. We then
obtain $$z(t+1/r)=M(x)\tau z(t)\tau^{-1}$$ for all $t\geq0$.
\par
For a twisted periodic continuously differentiable $g$ with $g(t)\in Gl_n(\C)$
for all $t\geq0$ let us denote
\begin{align*}
z_g(t) & = g(0)z(t)g^{-1}(t) \\
x_g(t) & = g(t)x(t)g^{-1}(t)-g'(t)g^{-1}(t).
\end{align*}
Then we have the following proposition.
\begin{Prop}\label{twistedDGL}
Let $x$ be a twisted periodic, continuous, matrix valued function, and let $z$
be the fundamental solution of $z'=zx$. Then
\renewcommand{\theenumi}{\roman{enumi}}
\renewcommand{\labelenumi}{\textnormal{(\theenumi)}}
\begin{enumerate}
\item
$z_g(t)$ is the fundamental solution of $z'_g=z_gx_g$.
\item
$M(x_g)=g(0)M(x)\tau g^{-1}(0)\tau^{-1}$.
\item
If $x_1$ is twisted periodic and there exists a $g_0$ such that
$$M(x_1)=g_0M(x)\tau g_0^{-1}\tau^{-1},$$
then there exists a twisted periodic Matrix $g(t)$ such that $g(0)=g_0$ and
$x_g(t)=x_1(t)$ for all $t\geq0$.
\end{enumerate}
\end{Prop}
\Proof
(i)\quad This is a direct calculation using  $(g^{-1})'=-g^{-1}g'g^{-1}$.
\par
(ii)\quad Since $g$ is twisted periodic, we have
$$M(x_g)=z_g(1/r)=g(0)z(1/r)g^{-1}(1/r)=g(0)M(x)\tau g^{-1}(0)\tau^{-1}.$$
\par
(iii)\quad Put $g(t)=z_1^{-1}(t)g_0z(t)$, where $z$ and $z_1$ are fundamental
solutions of $z'=zx$ and $z_1'=z_1x_1$ respectively. Then
the same calculation as in the proof of \cite{Frenkel}, Prop.(3.2.5) yields
$x_g = x_1$.
%
%
%
%
Again, a similar explicit calculation as in \cite{Frenkel} gives
$g(t+1/r) = \tau g(t) \tau^{-1}$ for all $t\geq0$.
\Endproof

Before we can use the results above to classify the affine adjoint orbits for
arbitrary $\g$  we need a  general fact from differential geometry, which is
stated as follows in
\cite{Frenkel}.

\begin{Prop}
Let $\g_\C\subset M_n(\C)$ be a matrix Lie algebra, and $G_\C\subset Gl_n(\C)$
the corresponding Lie group.  If $z$ is a solution of the linear differential
equation $z'=zx$, then $z(t)\in G_\C$ for all $t\geq0$ if and only if
$x(t)\in\g_\C$
for all $t\geq0$.
\end{Prop}

Now let $\tilde \g_\C$ be an affine Lie algebra of type $X_n^{(r)}$, and let
$\tau$ be the corresponding diagram automorphism of the underlying finite
dimensional Lie algebra $\g_\C$ used in the loop realization of $\tilde
\g_\C$.
In the case of an untwisted affine Lie algebra, $\tau$ is just the identity on
$\g_\C$. Let $\tilde \g$ and $\g$ denote the corresponding compact forms.
Following \cite{Frenkel},  we define an affine shell ("standard paraboloid" in
\cite{Frenkel}) to be the following  submanifold of codimension 2 in
$\tilde \g_\C$.
$$\P^{a,b}_\C=\{x(\cdot) + a_1C + b_1D\in\tilde \g_\C|2a_1b_1+(x,x)=a,
b_1=b\},$$
where $a,b\in\C$ and $b\neq0$. The zero-hyperplane in $\tilde \g_\C$ is
defined to
be the subspace $$\hat\g_\C=\{x(\cdot)+aC+bD\in\tilde\g_\C|b=0\}.$$
\par
By $\P^{a,b}$ with $a,b\in\R$, $a\neq0$ and $\hat \g$ we shall denote the
corresponding submanifolds of $\tilde\g$.  Let $\O_X$ be the $\L
G_{(\C)}$-Orbit
of $X$ in $\tilde\g_{(\C)}$, and let $\O_{g\tau}$ be the $G_{(\C)}$-Orbit of
$g\tau$ in $G_{(\C)}\tau$. Here $G\tau$ is the connected component of the
principal
extension $\tilde G$  of the compact group $G$ as constructed in
\S\ref{extension},
and $\tilde G_\C$ is the corresponding  complexification.

\begin{Theorem}\label{classification}
\renewcommand{\theenumi}{\roman{enumi}}
\renewcommand{\labelenumi}{\textnormal{(\theenumi)}}
\begin{enumerate}
\item
Each $\L(G_\C,\tau)$ (resp. $\L(G,\tau)$)-Orbit in the complex (resp. compact)
affine Lie algebra $\tilde\g_\C$ (resp. $\tilde\g$) lies either in one of the
affine shells $\P^{a,b}_\C$ (resp. $\P^{a,b}$) or in the zero-hyperplane.
\item
For a fixed affine shell, the monodromy map 
$$\O_{x(\cdot)+aC+bD}\mapsto\O_{M(\frac1bx)\tau}$$ 
is well defined and injective.
\item
For a fixed affine shell, the map defined in \textnormal{(ii)} gives a
bijection
between the $\L(G_\C,\tau)$ (resp. $\L(G,\tau)$)-Orbits in $\tilde\g_\C$
(resp. $\tilde \g$) which contain a constant loop and the $G_\C$ (resp.
$G$)-Orbits
in $G_\C\tau$ (resp. $G\tau$)  which contain an element that is invariant under
conjugation with $\tau$.
\end{enumerate}
\end{Theorem}
\Proof
(i) \quad Follows from the formula in lemma \ref{ad}. Note that the affine
shells
in a twisted affine Lie algebra  are just the  $\tau$-invariant parts of the
affine shells in the corresponding untwisted algebras.
\par
(ii)\quad
We look at the map $\O_{x(\cdot)+aC+bD}\mapsto \O_{M(\frac1bx)\tau}$ where,
as above,
$M(\frac1bx)=z(1/r)$, and $z$ is the fundamental solution of
$z'=z\cdot\frac1b\cdot x$.
Now, by (\ref{ad}),  we have 
$$\widetilde{\Ad}(g)(aC+bD+x)=\tilde aC+dD+gxg^{-1}-bg'g^{-1}.$$
But using \ref{twistedDGL} (i), we see that $\O_{\tilde
aC+bD+gxg^{-1}-bg'g^{-1}}$ is
being mapped to $\O_{z_g(1/r)\tau}$, and \ref{twistedDGL} (ii) yields
$$z_g(1/r)=M(\frac1bx_g)=g(0)M(\frac1bx)\tau g(0)^{-1}\tau^{-1},$$
hence 
$$z_g(1/r)\tau=g(0)z(1/r)\tau g(0)^{-1}\in\O_{z(1/r)\tau}.$$
So the map is well defined and injectivity follows with \ref{twistedDGL} (iii).
\par
(iii) \quad
If $s\tau\in G_{(\C)}\tau$ is invariant under conjugation with $\tau$ then
so is
$r\cdot b\cdot\log(s)$ and $\O_{a_1C+bD+r\cdot b\cdot\log(s)}$ is a preimage of
$\O_{s\tau}$ whenever it belongs to $\P^{a,b}$. On the other hand, if the
orbit
$\O_{a_2C+bD+x(\cdot)}$ contains a constant loop $aC+bD+x_0$, then $x_0$
has to
be invariant under conjugation with $\tau$. Now the fundamental solution of the
differential equation $z'=z\cdot \frac1bx_0$ is given by
$z(t)=\exp(t\frac1bx_0)$.
Hence $z(1/r)$ is invariant under conjugation with $\tau$ as well.
\Endproofshort

\begin{Cor}
If $\ord(\tau)=1$, or $G$ is compact, then the monodromy map in Theorem
\ref{classification} is surjective and hence induces a bijection between the
$\L G$-orbits in a fixed affine shell $\P^{a,b}$ and the $G$-orbits in $G\tau$.
\end{Cor}
\Proof
If $\tau=id$ then the statement is trivial. For compact $G$ we use the fact
that
every $G$-orbit in $G\tau$ intersects $S\tau$ in at least one point. Here
$S$ is
a Cartan subgroup of $\tilde G$ containing $\tau$.
\Endproofshort

\begin{Rem}
In the case of complex groups the classification of $\L (G_\C,\tau)$-orbits in
$\tilde \g_\C$ remains open. In this case it is no longer true that every
$G_\C$-orbit
in $G_\C\tau$ contains a $\tau$-invariant element (cf. \cite{Mohrdieck} for
an example),
so different arguments may have have to be applied. In any case, 
it is easy to see that 
the image under the monodromy map defined in \ref{classification} 
\textnormal{(ii)} are the 
$G_{(\C)}$-Orbits in $G_{(\C)}\tau$  for which there exists a $C^\infty$ path
$z:[0,1]\to G_{(\C)}$ such that $z(0)=e$, $z(1/r)\tau\in\O_{g\tau}$ and
$z(t+1/r)=z(1/r)\tau z(t)\tau^{-1}$ for all $t\geq0$.
\end{Rem}

\begin{Rem}
We have not dealt with the $\L(G,\tau)$-orbits in the zero-hyperplane, which
are basically  the orbits of the adjoint representation of $\L(G,\tau)$ 
on its Lie algebra. As we shall see later, the orbits relevant 
for representation theory are 
the $\L(G,\tau)$-orbits in some fixed affine shell with $b\neq 0$. Also, the 
classification of $\L(G,\tau)$-orbits in $\hat \g$ is presumably not 
manageable, i.e. it certainly yields an infinite dimensional ``moduli space''. 
\end{Rem}

\par
In a simply connected, semi-simple, compact Lie group, the conjugacy classes
are in one-to-one correspondence with a fundamental domain of the  affine Weyl
group $\tilde W$ acting on $LT$ in the notation of \S\ref{integration}. Here
the affine Weyl group acts like the finite Weyl group of  the root system $R$,
belonging to the group $G$, extended  by the group of translations generated by
the dual roots $R^*\in LT$. Here, the set of dual roots is given as
$R^*=\{\alpha^*\enspace |\enspace\alpha\in R\}$, and $(\alpha^*,.)=\alpha$
where
$(.,.)$ is the Killing form. Hence the the orbits in  a given affine shell
$\P^{a,b}$ are in one-to-one correspondence with a fundamental domain of
$\tilde W$.
\par
The $G$-orbits in $G\tau$ are in one-to-one correspondence with the set
$S_0\tau/W(S)$,
where $W(S)=N(S)/S$ acts by conjugation on $S_0\tau$. We can pull back this
action to
$S_0$ by right multiplication with $\tau^{-1}$. Then $tS_0\in W(S)$ acts on
$S_0$
via $s\mapsto ts\tau t^{-1}\tau^{-1}$.  By Proposition \ref{weylgruppen},
we have
$W(S)=(T/S_0)^\tau \rtimes W^\tau$. In this way, the action of $W^\tau$ is
the usual Weyl
group action  of the Weyl group of $G^\tau$ on the maximal torus $S_0$.
Hence the orbits
of this action are parametrized by a fundamental domain of the affine Weyl
group $\tilde
W^\tau$, where the translation part of $\tilde W^\tau$ is given by the dual
roots
$R^{\tau*}\in LS_0$ of the root system $R^\tau$ belonging to the group
$G^\tau$.
\par
An element $tS_0\in(T/S_0)^\tau$ acts on $S_0$ via $s\mapsto ts\tau
t^{-1}\tau^{-1}=
ss_{t^{-1}}\tau$, where $\tau t\tau^{-1}=ts_t$ with $s_t\in S_0$.
Viewing $LS_0$ as the universal covering of $S_0$ via $\exp$, we see that
$(T/S_0)^\tau$
acts on $LS_0$ as a group of translations.
A direct calculation shows that a set of generators of this group is given
by the set
$$\{\frac1r\sum_{i=1}^r\tau^i(\alpha^*)\enspace|\enspace\alpha^*\in
R^*,\tau(\alpha^*)\neq\alpha^*\}.$$
So together these two groups yield exactly the action of the affine Weyl group
$\tilde W^1$ belonging to the root system $R^1$ from \S\ref{integration}.
Hence we have proved
\begin{Prop}\label{weylgroupsII}
In the twisted compact case, the set of orbits in a given affine shell
$\P^{a,b}$ is
in one-to-one correspondence with a simplex in  $LS_0$ which is a
fundamental domain
for the action of the affine Weyl group $\tilde W^1$ belonging to the root
system $R^1$.
\end{Prop}

\section{Orbital integrals}\label{character}

\subsection{Affine Lie algebras and the Kac-Weyl 
character formula}\label{affcharacters}

Before we start deriving the analogue of Frenkel's character formula for the
twisted affine Lie algebras, let us briefly review some facts from the
structure
and representation  theory of affine Lie algebras. Let $\tilde\g_\C$ be an
affine Lie
algebra of type  $X_n^{(r)}$ with Weyl group $\tilde W$. Let $\tilde\h_\C$
be a Cartan
subalgebra and
$\Pi= \{\alpha_0,\ldots,\alpha_n\}\subset\tilde\h^*_\C$ be a set of simple
roots of
$\tilde\g_\C$ where the roots are labeled in the usual way (cf.\cite{Kac}).
$\tilde R$
denotes the set of roots of $\tilde\g_\C$, and define $R^\circ$ to be the
root system
which is obtained by deleting the 0-th vertex from the extended Dynkin
diagram of
$\tilde R$. Also, let $W^\circ$ denote the Weyl group  belonging to
$R^\circ.$ Let
$\alpha_i^\vee\in\tilde\h_\C$ be the dual simple roots such that
$\langle\alpha_i,\alpha_j^\vee\rangle=(A)_{i,j}$, where $A$  is the
generalized Cartan
matrix  belonging to $\tilde\g_\C$. Let $a_i$ be the "minimal" integers
such that
$A(a_0,\ldots, a_n)=0$, set $\delta=\sum_{i=0}^na_i\alpha_i$ and define
$d=\frac1{2\pi i}D$ with $D$ as in \S\ref{topologies}. Then we have
$\langle\alpha_i,d\rangle=0$ and $\langle\delta,d\rangle=1$.
Now, we can define an element $\theta=\delta-a_0\alpha_0\in \h_\C^{\circ*}$
and
the lattice $M= \Z W^\circ\theta^\vee\subset \h_\C^\circ$. Here
$\h_\C^\circ\subset\tilde\h_\C$ denotes the subspace generated by
$\{\alpha_1^\vee,\ldots\alpha_n^\vee\}$, and $\theta^\vee$ means the
element dual  to
$\theta$  in the sense of \cite{Kac}. A fundamental result in the theory of
affine Lie
algebras states that $\tilde W=W^\circ\ltimes M$. Observe that in
\cite{Kac} Kac uses
a lattice $M'\subset \h_\C^{\circ*}$ after identifying $ \h_\C^{\circ*}$ with
$ \h_\C^{\circ}$ via an invariant bilinear form.

\medskip

Turning to the representation theory of affine Lie algebras we define as usual
\begin{align*}
\tilde P=&\{\lambda\in\tilde \h^*_\C\enspace |\enspace
\langle\lambda,\alpha_i^\vee\rangle\in\Z\text{ for all }i=0,\ldots,n\}, \\
\tilde P_+=&\{\lambda\in\tilde P\enspace|\enspace
\langle\lambda,\alpha_i^\vee\rangle\geq0\text{ for all }i=0,\ldots,n\},
\text{ and }\\
\tilde P_{++}=&\{\lambda\in\tilde P\enspace|\enspace
\langle\lambda,\alpha_i^\vee\rangle>0\text{ for all }i=0,\ldots,n\}.
\end{align*}
Then there exists a bijection between the irreducible integrable highest
weight modules of
$\tilde\g_\C$ and the dominant integral weights $\Lambda\in\tilde P_+$. If
$L(\Lambda)$
is an integrable irreducible highest weight module with highest weight
$\Lambda\in\tilde
P_+$ then the formal character of $L(\Lambda)$ is the formal sum
$$\textnormal{ch}\enspace L(\Lambda)=\sum_{\lambda\in\tilde\h^*}\dim
L(\Lambda)_\lambda e(\lambda).$$
Here $ L(\Lambda)_\lambda$ is the weight space corresponding to the weight
$\lambda$,
and $e(\lambda)$ is a formal exponential.
The Kac-Weyl character formula now reads (cf. \cite{Kac})
$$\textnormal{ch}\enspace L(\Lambda)=\frac{\sum_{w\in\tilde W}
\epsilon(w)e(w(\Lambda+\tilde\rho))}
{e(\tilde\rho)\prod_{\alpha\in R}(1-e(-\alpha))^{mult(\alpha)}}.$$
Here $\tilde\rho$  is defined by
$\langle\tilde\rho,\alpha_i^\vee\rangle=1$ for $i=0,\ldots,n$ and
$\langle\tilde\rho,d\rangle=0$.
As usual we have set $\epsilon(w)=(-1)^{lenght(w)}$.
\par
So far, the character was considered as a formal sum involving the formal
exponentials
$e(\lambda)$.  Now we set $e^\lambda(h)=e^{\langle\lambda,h\rangle}$ for
$h\in\tilde\h_\C$.  In this way one can consider the character of a highest
weight
$\tilde\g_\C$-module $V$ as an infinite  series.  Let us set
$Y(V)=\{h\in\tilde\h_\C\enspace |\enspace \textnormal{ch}_V(h)\text{ converges
absolutely}\}$.   Then $ch_V$ defines a holomorphic function on $Y(V)$  and
the
following result holds \cite{Kac}:
\begin{Prop}
Let $V(\Lambda)$ be the irreducible highest weight module with highest
weight $\Lambda\in P_+$.
Then  $$Y(V(\Lambda))=\{h\in\tilde\h_\C\enspace |\enspace
Re\langle\delta,h\rangle>0\},$$
where $\delta=\sum_{i=0}^na_i\alpha_i$  and the $a_i$ are the labels of the
vertices of the
affine Dynkin diagram (cf. \cite{Kac}).
\par
In this setting the Kac-Weyl character formula gives an identity of
holomorphic functions
on $Y(V(\Lambda))$.
\end{Prop}
\par
Now we have $\tilde\h_\C=\h^\circ_\C\oplus\C C\oplus\C D$ with $C$
and $D$ as in \S\ref{topologies}.
With this notation one gets
$$Y(V(\Lambda))=\{h+aC+bD\enspace | \enspace h\in\h^\circ_\C,\enspace
a,b\in\C,\enspace Im\enspace b<0\}.$$

\subsection{Poisson transformation of the numerator of affine
characters}\label{poissontransformation}

In this section we will start do derive an analogue of Frenkel's character
formula for
twisted affine Lie algebras by deriving a formula for the numerator of the
character formula
in terms of the underlying non-connected Lie group.
To do this, we need to introduce some more notation. So let
$\tilde\g_\C$ be an arbitrary affine Lie algebra of type $X_n^{(r)}$, and
let $\tilde\g_\C'$ be
the untwisted affine Lie algebra of type  $X_n^{(1)}$ such that
$\tilde\g_\C\subset \tilde\g_\C'$ and let $\tilde\g$ and $\tilde\g'$ be the
corresponding
compact forms. If $\tilde\g_\C$ is untwisted, we have $\tilde\g_\C =
\tilde\g_\C'$. Furthermore,
let $R$ be the root system of the finite dimensional Lie algebra $\g_\C$
used to construct
$\L(\g_\C,\tau)$ in \S\ref{topologies}, and for a diagram automorphism
$\tau$ of $\g_\C$ let
$R^\tau$ be the "folded" root system introduced in
\S\ref{integrationformula}. If $\tilde\g_\C$
is a twisted affine Lie algebra with root system $\tilde R$, then
$R^\tau=R^\circ$. Also,
let $(.,.)$ denote the Killing form on $\tilde\g_\C'$, and let $(.,.)_r$
denote its restriction
to $\tilde\g_\C$.
\par
We now turn to the analytic Kac-Weyl character formula from
\S\ref{affcharacters}. Let
$\Lambda$ be a highest weight of $\tilde\g_\C$. There is no essential loss
in generality
by assuming $\langle\Lambda,d\rangle = 0$.
Then, after identifying $\tilde \h_\C\cong\tilde\h_\C^*$
via  $(.,.)_r$ we can choose $a\in\C$ and $H\in\h_\C^\circ$ such that
$\Lambda+\tilde\rho=aD+H$.
The condition $\Lambda\in\tilde P_+$ implies $a\in i\R$, $Im\enspace a<0$
and $H\in i\h^\circ$.
The numerator of the Kac-Weyl character formula  evaluated at $bD+K$ with
$b\in\C$, $Im\enspace b<0$
and $K\in\h_\C^\circ$ now reads
$$\sum_{w\in\tilde W}\epsilon(w)e^{(w(aD+H),bD + K)_r}.$$
Let us assume in the following $b\in i\R$ and $K\in i\h^\circ$. We then set
$t=-1/ab$, $h=H/a$
and $k=K/b$ yielding $t\in\R_+$ and $h,k\in\h$.
With $2\pi i d=D$,  the sum above reads
$$\sum_{w\in\tilde W}\epsilon(w)e^{-\frac 1t(w(2\pi i d+h),2\pi i d + k)_r}.$$
Let us set $c=2\pi i C$. Then the lattice $M$ operates on on $\tilde h_\C$ via
$$\gamma(h+ac+bd)=h+ac+bd-\left((h,\gamma)_r+\frac{ba_0}{2}\Vert\gamma\Vert_r^2\right)c.$$
Now a short calculation (cf. \cite{Frenkel}) shows for $w\in W^\circ$,
$\gamma\in M$, and
$w^{-1}\gamma\in\tilde W$ $$\left(w^{-1}\gamma(2\pi i d+h),2\pi i d +
k\right)_r =
-\frac12\Vert 2\pi i a_0\gamma+h-w k\Vert^2_r+\frac12\Vert
h\Vert^2_r+\frac12\Vert k\Vert^2_r.$$
Inserting this into the sum above, we get
\begin{multline*}\sum_{w\in\tilde W}\epsilon(w)e^{-\frac 1t(w(2\pi i
d+h),2\pi i d + k)_r} = \\
e^{-\frac{1}{2t}\Vert h\Vert^2_r-\frac{1}{2t}\Vert
k\Vert_r^2}\sum_{\gamma\in 2\pi i a_0 M}
\sum_{w\in  W^\circ}\epsilon(w)e^{\frac{1}{2t}\Vert\gamma+h-wk\Vert^2_r}.
\end{multline*}
We now need to apply the Poisson transformation formula: For a Euclidean
vector space $V$, a
lattice $Q\in V$ and a Schwartz function $f:V\to\C$ one has
$$\sum_{\mu\in Q^\vee}\hat f(\mu)=vol\enspace Q\sum_{\gamma\in Q}f(\gamma)$$
with
$$\hat f(\mu)=\int_Ve^{2\pi i(\gamma,\mu)}f(\gamma)d\gamma.$$
For a fixed $x\in \h^\circ$ set $f(\mu)=e^{(x,\mu)_r-\frac
t2\Vert\mu\Vert_r^2}$. Then we get
$$\hat f(\gamma)=\left(\frac{2\pi}{t}\right)^\frac l2e^{\frac1{2t}\Vert
x+2\pi i\gamma\Vert_r^2}.$$
with $l=\dim_\R\h^\circ$.
So for $x=h-wk$ with $h,k\in\h^\circ$ and $w\in W^\circ$ we obtain the identity
\begin{multline*}
\sum_{\gamma\in a_0 M}e^{\frac{1}{2t}\Vert 2\pi i\gamma+h-wk\Vert^2_r}=\\
 vol(a_0M)^{-1}\left(\frac{2\pi}{t}\right)^{-\frac l2} \sum_{\mu\in
(a_0M)^\vee}
e^{(\mu,h-wk)_r-\frac t2\Vert\mu\Vert^2_r}.
\end{multline*}
\par
If $\tilde R$ is of type Aff1, then $\theta$ is a long root in $R^\circ$,
and if $\tilde R$
is a root system of type Aff2 or Aff3, but not of type $\A_{2n}^{(2)}$,
then $\theta$ is a short
root in $R^\circ$. In case $\tilde R$ is of type $\A_{2n}^{(2)}$ then
$R^\circ$ is a root system
of type $\B\rC_n$, and $\theta$ is a root of medium length in $R^\circ$.
\par
So if $\tilde R$ is of type $X_n^{(r)}$ with $r=2,3$ and $\tilde R\neq
A_{2n}^{(2)}$, then
$\theta^\vee$ is a long root in $R^{\circ\vee}$, and hence $M$ is the
lattice which is generated
by the long roots in $R^{\circ\vee}$. If $\tilde R$ is of type
$A_{2n}^{(2)}$, then $\theta^\vee$
is a root of medium length in $R^{\circ\vee}$ and in this case we have
$a_0=2$. Thus, in all cases,
$M$ is the lattice which is generated by the  root system $R^1$ from
\S\ref{integrationformula}.
\par
For an arbitrary root system $S$ let $P^\circ(S)$ denote the weight lattice
of $S$. Then the above
implies
\begin{multline*}
\sum_{\gamma\in a_0 M}e^{\frac{1}{2t}\Vert 2\pi i\gamma+h-wk\Vert^2_r}=\\
\frac{1}{vol\enspace a_0M}\left(\frac{2\pi}{t}\right)^{-\frac l2}
\sum_{\mu\in P^\circ(R^1)}e^{(\mu,h-wk)_r-\frac t2\Vert\mu\Vert^2_r}.
\end{multline*}
Putting the above formulas together, we get
\begin{multline*}
\sum_{w\in\tilde W}\epsilon(w)e^{-\frac 1t(w(2\pi i d+h),2\pi i d + k)_r} = \\
\frac{e^{-\frac{1}{2t}\Vert h\Vert^2_r-\frac{1}{2t}\Vert k\Vert^2_r}}{vol(\Z\
R^1)\left(\frac{2\pi}{t}\right)^{\frac l2}}
\sum_{w\in W^\circ}\sum_{\mu\in P^\circ
(R^1)}\epsilon(w)e^{(\mu,h-wk)_r-\frac t2\Vert\mu\Vert^2_r}.
\end{multline*}
\par
Now let $W(R^1)$ denote the Weyl group of the root system $R^1$. It is a
well known fact that
(after the choice of a basis of $R^1$) every
weight $\lambda\in P^\circ(R^1)$ is conjugate under $W(R^1)$ to some
dominant weight
$\lambda'\in P_+^\circ(R^1)$. Since the root systems $R^1$
and $R^\circ$ are dual to each other, we have $W^\circ=W(R^1)$. So we get
\begin{multline*}
\sum_{w\in W^\circ}\sum_{\mu\in P^\circ
(R^1)}\epsilon(w)e^{(\mu,h-wk)_r-\frac t2\Vert\mu\Vert^2_r} = \\
\sum_{\mu\in P^\circ_+ (R^1)} \sum_{w\in W^\circ}\sum_{w'\in
W^\circ}\epsilon(ww')\epsilon(w')
e^{\langle w'\mu,h-wk\rangle}e^{-\frac t2\Vert\mu\Vert^2_r}
\end{multline*}
Here we have identified $\h^\circ$ and $\h^{\circ*}$ via $(.,.)_r$. In the
equation above, the
singular weights cancel out, so it is enough to
sum over the strictly dominant weights, or equivalently to replace
$\lambda$ by $\lambda+\rho^\tau$
with  $\rho^\tau=\frac12\sum_{\bar\alpha\in R^1_+}\bar\alpha$ as in
\S\ref{integrationformula}. Hence
\begin{align*}
\sum_{w\in W^\circ}&\sum_{\mu\in P^\circ
(R^1)}\epsilon(w)e^{(\mu,h-wk)_r-\frac t2\Vert\mu\Vert^2_r} = \\
& = \sum_{\lambda\in P^\circ_+ (R^1)} \sum_{w\in W^\circ}\sum_{w'\in
W^\circ}\epsilon(ww')\epsilon(w')
    e^{<w'(\lambda+\rho^\tau),h-wk>}e^{-\frac
t2\Vert\lambda+\rho^\tau\Vert^2_r}  \\
& =  \sum_{\lambda\in  P^\circ_+
(R^1)}\delta^\tau(h)\delta^\tau(-k)\chi_\lambda^\tau(h)\chi_\lambda^\tau(-k)
   e^{-\frac t2\Vert\lambda+\rho^\tau\Vert^2_r}
\end{align*}
with $A^\tau(\lambda)$, $\delta^\tau$ and $\chi_\lambda^\tau =
A^\tau(\lambda+\rho^\tau)/\delta^\tau$
as in \S\ref{applications}.
As before, let $\g_\C$ be the finite dimensional complex Lie algebra used
to construct $\L(\g_\C,\tau)$
with root system $R$ and compact form $\g$.
Let $G$ be the simply connected compact Lie group belonging to $\g$ and let
$G\tau$ denote the
connected component of the non-connected
Lie group $G\ltimes\langle\tau\rangle$ containing $\tau$. In
\S\ref{applications} we have seen
that $\chi_\lambda^\tau(h)=\tilde\chi_\lambda(e^h\tau)$
for $h\in\h^\circ$. Here $\tilde\chi_\lambda$ denotes the character of
$G\rtimes\langle\tau\rangle$
belonging to the highest weight $\lambda$
(cf. Theorem \ref{Charakterformel} and observe that in this notation we
have $\h^\circ=LS_0$ with
$LS_0$ as in
\S\ref{integrationformula}).
\par
So putting everything together, we have proved the following theorem (which
is the analogue of
Theorem (4.3.4) in \cite{Frenkel}).
\begin{Theorem}\label{poisson}
For $h,k\in\h^\circ$ one has
\begin{multline*}
e^{\frac{1}{2t}\Vert h\Vert^2_r}e^{\frac{1}{2t}\Vert
k\Vert^2_r}\sum_{w\in\tilde W}\epsilon(w)
e^{-\frac 1t(w(2\pi id+h),2\pi id+k)_r} =  \\
  \frac{\delta^\tau(h)\delta^\tau(-k)}{vol(\Z R^1)(\frac{2\pi}{t})^{\frac
l2}}\sum_{\lambda\in
P^\circ_+(R^1)}
   \chi_\lambda(e^h\tau)\chi_\lambda(e^{-k}\tau)e^{-\frac
t2\Vert\lambda+\rho^\tau\Vert^2_r}
\end{multline*}
\end{Theorem}

\begin{Rem}
This section is basically a reformulation of of the analogous results
of \cite{Frenkel} in the non-twisted cases  and of \cite{Kleinfeld} 
in the twisted cases. 
Kleinfeld  did his calculations for the twisted affine 
Lie algebras in concrete 
realizations of the corresponding root systems. 
Not realizing  the appearance of the characters of the non-connected Lie 
group $G\rtimes\langle\tau\rangle$, he had to work with 
the irreducible characters of
the different Lie groups corresponding to the root systems 
$R$, $R^\tau$ and $R^{\tau\vee}$.
\end{Rem}

\subsection{The heat equation}\label{heat}

In this section we will see how the expression for the numerator in Theorem
\ref{poisson} is connected
to the fundamental solution of the heat equation on the  component $G\tau$
of the non-connected group
$G\rtimes\langle\tau\rangle$. To this end, let $\Delta_G$ denote the
Laplacian on the compact simply
connected group $G$ with respect to the Riemannian metric on $G$ induced by
the negative of the Killing
form on $\g_\C$. We can pull back this metric to $G\tau$ such that right
multiplication with $\tau$
induces an isometry between the Riemannian manifolds $G$ and $G\tau$. The
Laplacian on $G\tau$ shall
be denoted with $\Delta_{G\tau}$.
\par
Now for a fixed parameter $T>0$, the heat equation on $G\tau$ reads
$$\frac{\partial f(g\tau,t)}{\partial t}=\frac {sT}2\Delta_{G\tau} f(g\tau,t)$$
with $g\in G$, $s\in\R$, $s>0$ and where $f:G\tau\to\R$ is continuous in
both variables, $C^2$ in the
first and $C^1$ in the second variable.
The fundamental solution of the heat equation is defined by the initial data
$$ f(g\tau,t)|_{t=+0}=\delta_\tau(g\tau),$$
where $\delta_\tau$ is the Dirac delta distribution centered at $\tau\in
G\tau$.
\par
For a highest weight $\lambda\in P^\circ_+(R)$ of $G$ let $d(\lambda)$
denote the dimension of the
corresponding irreducible representation of $G$ and
$\chi_\lambda$ its character. Then the fundamental solution of the heat
equation on $G$ is given
by  
$$u_s(g,t)=\sum_{\lambda\in
P^\circ_+(R)}d(\lambda)x_\lambda(g)
e^{-\frac{stT}2(\Vert\lambda+\rho\Vert^2-\Vert\rho\Vert^2)}$$
(see \cite{Fegan}). Now $G$ and $G\tau$ are isometric as Riemannian 
manifolds, hence the
fundamental solutions of the
respective heat equations coincide. That is, the fundamental solution of
the heat equation on
$G\tau$ is given by
$$v_s(g\tau,t)=\sum_{\lambda\in P^\circ_+(R)}d(\lambda)x_\lambda(g)
e^{-\frac{stT}2(\Vert\lambda+\rho\Vert^2-\Vert\rho\Vert^2)}.$$
\par
There is a well known identity for the characters of a compact group which
can be derived using the
orthogonality relations for irreducible characters:
$$d(\lambda)\int_G\chi_\lambda(g_1gg_2^{-1}g^{-1})dg=
\chi_\lambda(g_1)\chi_\lambda(g_2^{-1}),$$
where $dg$ denotes the normalized Haar measure on $G$. Using a version of
the orthogonality
relations for non-connected groups,
it is easy to proof an analogous formula for the characters on the outer
components:
$$d(\lambda)\int_G\chi_\lambda(g_1\tau
g\tau^{-1}g_2^{-1}g^{-1})dg=
\chi_\lambda(g_1\tau)\chi_\lambda(\tau^{-1}g_2^{-1}).$$
Hence we obtain
\begin{multline*}
\int_G v_s(gg_1\tau g^{-1}\tau^{-1}g_2^{-1}\tau,t)dg = \\
\sum_{\lambda\in
P^\circ_+(R)}\chi_\lambda(g_1\tau)\chi_\lambda(\tau^{-1}g_2^{-1})
e^{-\frac{stT}2(\Vert\lambda+\rho\Vert^2-\Vert\rho\Vert^2)}.
\end{multline*}
By Theorem \ref{Charakterformel} we see that $\chi_\lambda(g\tau)=0$ if
$\lambda$ is not
$\tau$-invariant.
Furthermore, we have
$\Vert\lambda+\rho\Vert^2=\Vert\lambda+\rho^\tau\Vert^2_r$ if $\lambda$ is
$\tau$-invariant.
Thus, using Theorem \ref{poisson}, exchanging the role of $s$ and $t$, and
fixing the parameter value
$s = T$, we have proved the following proposition.
\begin{Prop}\label{heatprop}
For $h,k\in\h^\circ$ one has
\begin{multline*}
\sum_{w\in\tilde W}\epsilon(w)e^{-\frac 1t(w(2\pi id+h),2\pi id+k)_r}  = \\
\frac{e^{-\frac{1}{2t}\Vert h\Vert^2_r}e^{-\frac{1}{2t}\Vert
k\Vert^2_r}e^{-\frac{t}{2}\Vert
\rho^\tau\Vert^2_r}\delta^\tau(h)\delta^\tau(-k)}
{vol(\Z R^1)(\frac{2\pi}{t})^{\frac l2}}\cdot \\
\int_{G}v_{\frac t{T^2}}(ge^h\tau
g^{-1}\tau^{-1}e^{-k}\tau,T)dg
\end{multline*}
\end{Prop}

\subsection{Wiener measures and a path integral}\label{wienermeasures}

The main result of this subsection will be a further reformulation of the
numerator of the character formula as an integral over a certain path space
on the connected component $G\tau$. This is based on Proposition
\ref{heatprop},
above, and the theory of Wiener measure on $G\tau$ which we will study first. 
(Compare with \cite{Frenkel} and consult e.g. \cite{Kuo} for a 
comprehensive treatment of the theory
of Wiener measures on a vectorspace.)
\par
\medskip
The Wiener measure on an euclidian vectorspace $V$ 
of variance $s>0$ is a measure $\omega_V^s$ on the Banach space of paths
$$C_V=\{x:[0,T]\to V\enspace | \enspace x(0)=0, \enspace x 
\text{ continuous}\}$$
(with the supremum norm) and is defined using the fundamental solution 
$w_s(x,t)$ of the heat equation 
$$\frac{\partial f(x,t)}{\partial t}=\frac {sT}2\Delta_{V} f(x,t).$$
on $V$ as follows:
First, one defines cylinder sets in $C_V$ to be the following subsets 
of $C_V$:
$$\{x\in C_V : (x(t_1)\in A_1,\cdots, x(t_m)\in A_m)\},$$
with $0<t_1\leq t_2\leq\ldots\leq t_m\leq T$, $m\in\N$, and where
$A_1,\ldots A_m$
are Borel sets in $V$. 
Then the Wiener measure $\omega^s_V$ of variance $s>0$ 
is defined on the cylinder sets
of $C_V$ via
\begin{multline*}
\omega^s_V(x(t_1)\in A_1,\cdots, x(t_m)\in A_m)=  \\
\int_{A_1}\cdots\int_{A_m}w_s(\Delta x_1,\Delta t_1)\cdots
w_s(\Delta x_m,\Delta t_m)
d x_1\cdots d x_m,
\end{multline*}
where $dx$ is a Lebesgue measure on $V$ and we have set
$x_k=x(t_k)$, $\Delta x_k=x_k-x_{k-1}$, 
$\Delta t_k=t_k-t_{k-1}$ and $x_0=0$.

The conditional Wiener measure $\omega_{V,Z}^s$ of variance
$s>0$ is defined on the closed subspace $C_{V,X}\subset C_V$ 
with fixed endpoint $x(T)=X$ on the 
cylinder sets via
\begin{multline*}
\omega^s_{V,X}(x(t_1)\in A_1,\cdots, x(t_{m-1})\in A_{m-1})=  \\
\int_{A_1}\cdots\int_{A_{m-1}}w_s(\Delta xt_1,\Delta t_1)\cdots
w_s(\Delta x_m,\Delta t_m)
d x_1\cdots d x_{m-1},
\end{multline*}
where additionally $x_m=X$ and $t_m=T$.
\par
Now a classical result in the theory of Wiener measures states that 
the measures  $\omega^s_V$ and  $\omega_{V,X}^s$ are $\sigma$-additive on the 
$\sigma$-field generated by the cylinder sets in $C_V$ and $C_{V,X}$ 
respectively. 
Furthermore, the $\sigma$-fields generated by the cylinder sets are 
exactly the Borel fields of the respective Banach spaces (cf. \cite{Kuo}).  
As another result, we have
$$\omega_V^s(C_V)=1$$
and
$$\omega_{V,X}^s(C_{V,X})=w_s(X,T),$$
where $w_s(x,t)$ is the fundamental solution of the heat equation on $V$.
\par
Using the fundamental solution $u_s(g,t)$ of the heat 
equation on the compact Lie group $G$, 
one can define the Wiener measure $\omega_G^s$ and the conditional 
Wiener measure 
$\omega_{G,Z}^s$ on the complete metric spaces
$$C_G=\{z:[0,T]\to G\enspace | \enspace z(0) = e, 
\enspace z \text{ continuous}\}$$
and $C_{G,Z}=\{z\in C_G,\enspace z(T)=Z)\}$ in exactly the 
same fashion as the Wiener measure
on the vectorspace $V$ (cf. \cite{Frenkel}). 
The metric $\varrho$ on $C_G$ is given by 
$\varrho(z,z_1)=\sup_{t\in[0,T]}\varrho_0(z(t),z_1(t))$, 
where $\varrho_0(g,g_1)$ denotes the 
length of a shortest geo\-desic in $G$ connecting 
two given points $g$ and $g_1$ 
(the metric on $G$ still being given by the negative of 
the Killing form on $\g$).
\par
There is an important connection between the Wiener measure on $G$ and the 
Wiener measure on the Lie algebra $\g$ which was discovered 
by Ito \cite{Ito} and
explicitly constructed  by McKean \cite{McKean}:  
Let $y\in C_\g$ be a continuous path. 
Then for $n\in\N$ and $k=0,\ldots,2^n-1$, 
we define a path $z_n:[0,T]\to G$ by $z_n(0)=e$ and
$$z_n(t)=z_n(\frac k{2^n}T)\exp\left(y(t)-y(\frac k{2^n}T)\right)
\quad\text{for }\frac k{2^n}T< t\leq \frac{k+1}{2^n}T.$$
Note that if $y$ is a differentiable path, 
then $\lim_{n\to\infty}z_n$ is the fundamental 
solution of the differential equation of $z'=zy'$ 
and hence a well defined path in $G$. 
Let us define a map
$$i: C_\g\to C_G$$
via
$$y\mapsto
\begin{cases}
\lim_{n\to\infty} z_n\quad & \text{ if the limit exists,} \\
e  & \text{ else. }
\end{cases}$$
The fundamental result of Ito and McKean states that the series
$z_n$, converges with $n\to\infty$ in the topology of $C_G$ almost everywhere with 
respect to the measure $\omega_\g^s$.
Furthermore, the measure on $C_G$ induced by the map $i$ coincides with the Wiener 
measure $\omega_G^s$ on $G$. 
Hence $i$ induces an isomorphism 
$$I:L_1(C_\g,\omega_\g^s)\to L_1(C_G,\omega_G^s)$$
via $If(iy)=f(y)$. 
\par
This isomorphism is called Ito's isomorphism in the literature. With its
help it is easy to translate most of
the results about the Wiener measure on a vectorspace to a 
corresponding result about Wiener measure on a compact Lie group. 
For example we have
$$\omega_G^s(C_G)=1$$
and
$$\omega_{G,Z}^s(C_{G,Z})=u_s(Z,T),$$
where $u_s(g,t)$ denotes the fundamental solution of the heat equation on $G$ 
(cf. \cite{Frenkel}).
\par
We will denote the integrals with respect to the Wiener measure 
and the conditional Wiener measure  on $V$ as
$$\int_{C_V}f(x)d\omega^s_{V}(x),$$
resp.
$$\int_{C_{V,X}}f(x)d\omega^s_{V,X}(x),$$
and accordingly for the Wiener measures on $G$.
\par
One of the most important properties of these integrals is its
translation qasi-invariance (cf. \cite{Kuo}):
Let $f:C_V\to \R$ be an integrable function 
and let $y\in C_V$ be a $C^\infty$-path.
Then
$$\int_{C_V}f(x)d\omega_V^s(x) = 
\int_{C_V}f(x+y)
e^{-\frac 1s(x',y')-\frac 1{2s}(y',y')}d\omega_V^s(x).$$
For a function $f:C_{V,X}\to \R$ the  translation quasi-invariance of 
the conditional Wiener measure reads
$$\int_{C_{V,X}}f(x)d\omega_{V,X}^s(x) = 
\int_{C_{V,X+Y}}f(x+y)
e^{-\frac 1s(x',y')-\frac 1{2s}(y',y')}d\omega_{V,X+Y}^s(x),$$
with $Y=y(T)$. In the above formulas, $(x',y')$ denotes the 
Stieltjes integral $\frac 1T\int_0^T(y'(t),dx(t))$, and $(.,.)$ denotes 
the scalar product on $V$.
\par
Translated to the Wiener integral on $G$, 
the translation quasi-in\-var\-i\-ance 
looks as follows (see \cite{Frenkel}):
\begin{Prop}
\renewcommand{\theenumi}{\roman{enumi}}
\renewcommand{\labelenumi}{\emph{(\theenumi)}}
\begin{enumerate}
\item
Let $f:C_G\to\R$ be an integrable function, and $g\in C_G$ 
be a $C^\infty$-path.
Then
\begin{multline*}
\int_{C_G}f(z)d\omega_G^s(z) = \\
\int_{C_G}f(zg) e^{-\frac 1s(z^{-1}z',g'g^{-1})_\g - 
\frac 1{2s}(g^{-1}g',g^{-1}g')_\g}d\omega_G^s(z).
\end{multline*}
\item
Let $f:C_{G,Z}\to\R$ be an integrable function, 
and $g\in C_G$ be a $C^\infty$-path.
Then
\begin{multline*}
\int_{C_{G,Z}}f(z)d\omega_{G,Z}^s(z) = \\
\int_{C_{G,Zg(T)^{-1}}}f(zg)
e^{-\frac 1s(z^{-1}z',g'g^{-1})_\g}\cdot \\
e^{-\frac 1{2s}(g^{-1}g',g^{-1}g')_\g} d\omega_{G,Zg(T)^{-1}}^s(z).
\end{multline*}
\end{enumerate}
\end{Prop}
Here the term $(z^{-1}z',g'g^{-1})_\g$ should be 
interpreted as the Stieltjes integral
$\frac 1T\int_0^T(g'g^{-1},d(i^{-1}(z))_\g$, where $(.,.)_\g$ denotes the 
negative of the Killing form on $\g$ (we add the subscript $\g$ here to avoid
possible confusions in later calculations). Note that $i^{-1}$ 
is a well defined map almost everywhere on $C_G$ with respect to $\omega^s_G$.
\par
Using the translation quasi-invariance of the 
Wiener measure on $G$, Frenkel computes
the following integral with respect to this measure 
(\cite{Frenkel}, Prop. 5.2.12):
\begin{Lemma}\label{berechnung}
Let $Y\in \L(\g,\tau)$ and let $g\in C_G$ be a 
$C^\infty$- path such that $g'=gY$. Then
$$e^{-\frac{\Vert Y\Vert_\g^2}{2s}}
\int_{C_{G,Z}}e^{\frac 1s(z^{-1}z',Y)_\g}d\omega^s_{G,Z}(z) = 
u_s(Zg(T)^{-1},T),$$
where $u_s(z,t)$ is the fundamental solution of the heat equation on $G$.
\end{Lemma}
%
%
%
%
%
%

\par
Now let $\O_{g\tau}$ denote the $G$-orbit in $G\tau$ containing the 
element $g\tau$ (cf. \S \ref{orbitclassification}). Multiplying 
each element of $G\tau$ with $\tau^{-1}$, we can identify $\O_{g\tau}$ 
with a $G$-orbit in $G$, where $G$ acts on itself by twisted conjugation:
$(h,g)\mapsto hg\tau h^{-1}\tau^{-1}$. This orbit will be denoted with 
$\O_{g\tau}$ as well. We can now define $C_{G,\O_{g\tau}}\subset C_G$ to
be the space of continuous paths with $z:[0,T]\to G$ with $z(T)\in\O_{g\tau}$.
A conditional Wiener measure $\omega^s_{G,\O_{g\tau}}$  on  $C_{G,\O_{g\tau}}$
is defined via
\begin{multline*}
\int_{C_{G,\O_{g\tau}}}f(z)d\omega^s_{G,\O_{g\tau}}(z)= \\
  \int_{G}\left(\int_{C_{G,g_1g\tau g_1^{-1}\tau^{-1}}}f(z)
d\omega^s_{G,g_1g\tau g_1^{-1}\tau^{-1}}(z)\right)dg_1,
\end{multline*}
where $f$ is integrable on $C_{G,g_1g\tau g_1^{-1}\tau^{-1}}$ 
for almost all $g_1\in G$.
Inserting this definition into Lemma \ref{berechnung}, one gets:
\begin{Cor}\label{integCor}
Let $Y\in \L(\g,\tau)$ and let $g\in C_G$ be a 
$C^\infty$- path such that $g'=gY$. Then
\begin{multline*}
e^{-\frac{\Vert Y\Vert_\g^2}{2s}}
\int_{C_{G,\O_{Z\tau}}}e^{\frac 1s(z^{-1}z',Y)_\g}
d\omega^s_{G,\O_{Z\tau}}(z) = \\
\int_G u_s(g_1Z\tau g_1^{-1}\tau^{-1}g(T)^{-1},T)dg_1,
\end{multline*}
where $u_s(z,t)$ is the fundamental solution of the heat equation on $G$.
\end{Cor}

\par
We want to interpret the numerator of the Kac-Weyl character formula 
as an integral in a space of paths in $G\tau$, so we need to define 
a Wiener measure on the space 
$$C_{G\tau}=\{\tilde z:[0,T]\to G\tau \enspace | 
\enspace \tilde z\tau^{-1}\in C_G\}.$$
This can be done with the help of the fundamental solution $v_s(g\tau,t)$
on $G\tau$ in exactly the same way as on $G$ and on $V$. But as we have
seen before, we have $v_s(g\tau,t)=u_s(g,t)$ for all $g\in G$ and $t>0$.
So the integrals on $G$ and $G\tau$ will not differ, and we can define 
the measure $\omega^s_{G\tau}$ directly by 
$$\int_{C_{G\tau}}f(\tilde z)d\omega^s_{G\tau}(\tilde z) =  
\int_{C_{G}}\hat f(\tilde z\tau^{-1})d\omega^s_{G}(\tilde z\tau^{-1}),$$
 where $\hat f$ is a function on $C_G$ which is given by
$\hat f(z)=f(z\tau)$. The conditional Wiener measures on $C_{G\tau,Z\tau}$ 
and $C_{G\tau,\O_{g\tau}}$ are defined analogously.
\par
So the formula in Corollary \ref{integCor} now reads
\begin{multline*}
e^{-\frac{\Vert Y\Vert_\g^2T^2}{2t}}
\int_{C_{G\tau,\O_{Z\tau}}}e^{\frac {T^2}t (z^{-1}z',Y)_\g}
d\omega^\frac{t}{T^2} _{G\tau,\O_{Z\tau}}(z\tau) = \\
\int_G v_\frac{t}{T^2}(g_1Z\tau g_1^{-1}\tau^{-1}g(T)^{-1}\tau,T)dg_1,
\end{multline*}
where $v_s(g\tau,t)$ is the fundamental solution of the heat equation 
on $G\tau$.
Note, that we also have set $s=\frac{t}{T^2}$ in the above calculation.
\par
Now we fix the parameter value $T=\frac 1r$. Observe that for $Y\in\L(\g,\tau)$
we then have  $\Vert Y\Vert_r = -\Vert Y\Vert_\g$.
So for $h,k\in\h^\circ$, we can set $Y=\frac 1T k = rk$ and $Z=e^h$. Then 
we have $g(T)^{-1}=e^{-k}$. 
Hence the integral formula above and Proposition
\ref{heatprop} yield:
\begin{Prop}\label{wienercor}
Let $h,k\in \h^\circ$.
Then
\begin{multline*}
\sum_{w\in\tilde W}\epsilon(w)e^{-\frac 1t(w(2\pi id+h),2\pi id+k)_r} =  \\
\frac{\delta^\tau(h)\delta^\tau(-k)e^{-\frac{1}{2t}\Vert h\Vert_r-
\frac t2 \Vert\rho^\tau\Vert_r}}
{vol(\Z R^1)\left(\frac{2\pi}{t}\right)^\frac l2}
\int_{C_{G\tau,\O_{e^h\tau}}}e^{\frac 1{tr^2}(z^{-1}z',k)_\g}
d\omega^{tr^2}_{G\tau,\O_{e^h\tau}}(z\tau)
\end{multline*}
\end{Prop}

\subsection{Affine characters and orbital integrals}\label{interpretation}

In this section we will indicate, how the integral 
in Proposition \ref{wienercor}
can be interpreted as an integral over an affine coadjoint orbit of
$\L(G,\tau)$. This interpretation and the analytic Kac-Weyl character
formula from \S\ref{affcharacters} will then yield an analogue of the 
Kirillov character formula for compact semisimple Lie groups. For precise
details we have to refer to Frenkel's work \cite{Frenkel}.
\par
\medskip
For fixed $a,b\in \R$ and $b\neq0$ let $\P^{a,b}$ be the affine shell defined
in \S\ref{orbitclassification}. We can identify $\P^{a,b}$ with
$\L(\g,\tau)$ via the projection $p:C\mapsto 0$ and $D\mapsto 0$.
Under this projection, the affine adjoint $\L(G,\tau)$-action on
$\L(\g,\tau)$ is given by $(g,y)\mapsto gyg^{-1} - bg'g^{-1}$ 
(cf. Prop. \ref{ad}). 
\par
We have a series of maps 
$$\P^{a,b}\overset{p}\longrightarrow \L(\g,\tau)
\overset{s}{\longrightarrow} C_\g^\infty 
\overset{i}{\longrightarrow}
C_G^\infty \overset{e_\tau}{\to}G\tau,$$
with $s(x)(t)=\int_0^tx(\kappa)d\kappa,$ and where $i$ 
maps a path $y\in C^\infty_\g$
to the fundamental solution of the differential equation $z'=\frac1bzy'$.
The map $e_\tau$ is given by $e_\tau(z) = z(1/r)\tau$.
From Ito's isomorphism we have a map $\tilde i:C_\g\to C_G$, 
which is the extension of the map $i: C_\g^\infty \to 
C_G^\infty$ above to the corresponding completions $C_\g$ and $C_G$.
\par
Now every element $y\in C_\g$ defines an element $dy\in\L(\g,\tau)^*$ 
via the Stieltjes integral 
$$\langle x,dy\rangle = r\int_0^{\frac 1r}(x(\kappa),dy(\kappa)),$$
where $(.,.)$ denotes the Killing form on $\g$ as in \S\ref{topologies}
and \S\ref{affcharacters}.
Note that for $y\in\L(\g,\tau)$, we have 
$\langle x,d(s(y))\rangle=(x,y)_r$, where 
$(.,.)_r$ is the bilinear form on $\L(\g,\tau)$ defined in 
\S\ref{affcharacters}. Let $\L(\g,\tau)^*_0$ denote the image of $C_\g$
under this map with the topology induced from $C_\g$,
and let $\tilde s:\L(\g,\tau)^*_0\to C_\g$ denote the
inverse map. In this notation, $\tilde s$ is 
the extension of the map $s$ to $\L(\g,\tau)^*_0$. 
Putting the above remarks together, 
we get the following commutative diagram:
\begin{center}
$\begin{array}{cccccccc}
\P^{a,b}\overset{p}{\longrightarrow} & 
\L(\g,\tau) & 
\overset s\longrightarrow & 
C^\infty_\g & 
\overset{i}{\longrightarrow} & 
C^\infty_{G} & 
\overset{e_\tau}{\longrightarrow} &
G\tau   \\
& \bigcap & & \bigcap & & \bigcap &&\\
& \L(\g,\tau)_0^* & 
\overset{\tilde s}{\longrightarrow} & 
C_\g &
\overset{\tilde i}{\longrightarrow} & 
C_{G} & 
\overset{\tilde e_\tau}{\longrightarrow} &
G\tau   
\end{array}$
\end{center}
\par
Here $\tilde e_\tau$ denotes the extension of $e_\tau$ to $C_G$.
We have seen in \S\ref{orbitclassification} that under the 
composition $e_\tau\circ i\circ s\circ p$,
a $\L(G,\tau)$-orbit $\O_{x(\cdot)+a_1C+b_1D}\subset \P^{a,b}$ 
is mapped to the $G$-orbit $\O_{e_\tau\circ i\circ s(x)}\subset G\tau$. 
\par
\medskip
On the Hilbert space $\L(\g,\tau)(L_2)$ introduced in 
\S\ref{topologies}, we can define a norm by 
$$|x(\cdot)|=\sup_{t\in[0,\frac1r]}
\vert\int_0^tx(\kappa)d\kappa\vert.$$
The completion of $\L(\g,\tau)(L_2)$ with respect to this norm will be 
$\L(\g,\tau)_0^*$. So we have a series of completions
$$\L(\g,\tau)\subset\L(\g,\tau)(L_2)\subset\L(\g,\tau)_0^*,$$
with respect to the $L_2$-topology on  $\L(\g,\tau)$ and the
norm on $\L(\g,\tau)(L_2)$ introduced above.
\par 
In this picture, the set 
$\tilde s^{-1}\circ\tilde i^{-1}\circ\tilde e_\tau^{-1}
(\O_{e_\tau\circ i\circ s(x)})
\subset \L(\g,\tau)_0^*$ can be viewed as the closure of the 
affine adjoint orbit $\O_{x(\cdot)+a_1C+b_1D}$ in $\L(\g,\tau)_0^*$
and is mapped to $C_{G,\O_{i\circ s(x)}}$ 
under the map $\tilde i\circ\tilde s$. 
Accordingly, the integral 
$$\int_{C_{G\tau,\O_{Z\tau}}}f(z)
d\omega^{s}_{G\tau,\O_{Z\tau}}(z\tau)$$
can be viewed as an integral over the closure of the corresponding
affine adjoint orbit in  $\L(\g,\tau)_0^*$. For more details, which involve
the construction of a Gaussian measure on  $\L(\g,\tau)_0^*$, 
cf. \cite{Frenkel}
\par
The discussion above allows us to interpret the integral 
appearing in the formula for the numerator of
the Kac-Weyl character formula in Prop. \ref{wienercor}
as an integral over the closure in $\L(\g,\tau)_0^*$ of the affine 
adjoint orbit containing $aD+H=\Lambda+\tilde\rho$. 
The denominator of the Kac-Weyl character formula is a function $p$
not depending on the highest weight $\Lambda$ of the corresponding
representation and hence can be seen as the analogue of the 
universal function  appearing in the Kirillov character 
formula for compact Lie groups
(in the context of compact Lie groups, the universal 
function is given by the denominator of 
the Weyl character formula as well).
Hence our character formula for affine Lie algebras can be
written in the following way:
%
%
\begin{Theorem}\label{Kirillovcharformel}
Let $bD+K\in\tilde h$ with $K\in\h^\circ$ and $b\in i\R$, $im (b)<0$. 
Furthermore, for $\Lambda\in\tilde P_+$ let 
$\Lambda+\tilde\rho = aD+H$. Then the character of the highest 
weight representation corresponding to $\Lambda$ evaluated at
$bD+K$ is given by
\begin{multline*}
\textnormal{ch}\enspace L(\Lambda)(bD+K)=p^{-1}(bD+K)\cdot \\
\frac{\delta^\tau(\frac Ha)\delta^\tau(-\frac Kb)
e^{\frac{ab}{2}\Vert \frac Ha\Vert_r
-\frac{1}{2ab} \Vert\rho^\tau\Vert_r}}
{vol(\Z R^1)\left(-2ab\pi\right)^\frac l2}\cdot\\
\int_{C_{G\tau,\O_{(e^{\frac Ha})\tau}}}e^{-\frac{ab}{r^2}
(z^{-1}z',\frac Kb)_\g}
d\omega^{-\frac {r^2}{ab}}_{G\tau,\O_{(e^{\frac Ha})\tau}}(z\tau).
\end{multline*}
\end{Theorem}
Following \cite{Frenkel}, the path integral above can be interpreted as
a Gaussian integral over the closure of the affine orbit containing
$aD+H$ in $\L(\g,\tau)_0^*$. In that setup, the above formula may be seen
as an exact analogue of Kirillov's classical character formula.

\end{document}